\documentclass[a4paper,9pt]{article}
\usepackage[latin1]{inputenc}
\usepackage{graphicx}

\usepackage{amsmath}
\usepackage{amsthm}
\usepackage{amssymb}
\usepackage{color}
\usepackage{mathrsfs}
\usepackage{fancyhdr}
\usepackage{setspace}
\usepackage[center]{subfigure}
\usepackage{float}
\usepackage[center,small,bf,ruled]{caption}

\begin{document}

\numberwithin{equation}{section}

\newtheorem {teo} {Theorem} [section]
\newtheorem {cor} {Corollary} [section]
\newtheorem {lem}[teo] {Lemma}

\theoremstyle{definition}
\newtheorem {oss} {Remark} [section]

\newcommand {\Dim} {\textsc{Proof\\ }}
\newcommand {\Fine} { $\blacksquare$ \\}



\setcounter{page}{1}

\title{{\bf Travelling Randomly on the Poincar\'e Half-Plane with a Pythagorean Compass}}
\maketitle
\author{
{\center
{ \large V. Cammarota \;\;\;\;\;  E. Orsingher \footnote{Corresponding author. Tel.: +390649910585, fax: +39064959241.\\ {\it E-mail address}: enzo.orsingher@uniroma1.it . }} \\
\vspace{5mm}
{\small Dipartimento di Statistica, Probabilit\`a e Statistiche applicate \\ University of Rome `La Sapienza'  \\ P.le Aldo Moro 5, 00185 Rome, Italy
}

}

\begin{abstract}
A random motion on the Poincar\'e half-plane is studied. A particle runs on the geodesic lines changing direction at Poisson-paced times. The hyperbolic distance is analyzed, also in the case where returns to the starting point are admitted. The main results concern the mean hyperbolic distance (and also the conditional mean distance) in all versions of the motion envisaged. Also an analogous motion on orthogonal circles of the sphere is examined and the evolution of the mean distance from the starting point is investigated.
\end{abstract}


AMS Classification 60K99\\

\section{Introduction}

\hspace{5mm} Motions on hyperbolic spaces have been studied since the end of the Fifties and most of the papers devoted to them deal with the so-called hyperbolic Brownian motion (see, e.g., Gertsenshtein and Vasiliev \cite{get}, Gruet \cite{gruet}, Monthus and Texier \cite{mon}, Lao and Orsingher \cite{l}).

 More recently also works concerning two-dimensional random motions at finite velocity on planar hyperbolic spaces have been introduced and analyzed (Orsingher and De Gregorio \cite{DO}).

While in \cite{DO} the components of motion are supposed to be independent, we present  here a planar random motion with interacting components. Its counterpart on the unit sphere is also examined and discussed.

The space on which our motion develops is the Poincar\'e upper half-plane $H_2^+=\{(x,y): y>0\}$ which is certainly the most popular model of the Lobachevsky hyperbolic space. In the space $H^+_2$ the distance between points is measured by means of the metric
\begin{equation}
\mathrm{d}s^2=\frac{\mathrm{d}x^2+\mathrm{d}y^2}{y^2}. \label{metric} \end{equation}

The propagation of light in a planar non-homogeneous medium, according to the Fermat principle, must obey the law
\begin{equation}
\frac{\sin \alpha(y)}{c(x,y)}=\mathrm{cost}
\end{equation}
where $\alpha(y)$ is the angle formed by the tangent to the curve of propagation with the vertical at the point with ordinate $y$. In the case where the velocity $c(x,y)=y$ is independent from the direction, the light propagates on half-circles as in $H_2^+$.

In \cite{optics} it is shown that the light propagates in a non-homogeneous half-plane $H_2^+$ with refracting index $n(x,y)=1/y$ with rays having the structure of half-circles.

Scattered obstacles in the non-homogeneous medium cause random deviations in the propagation of light and this leads to the random model analyzed  below.

The position of points in $H_2^+$ can be given either in terms of Cartesian coordinates $(x,y)$ or by means of the hyperbolic coordinates $(\eta, \alpha)$. In particular, $\eta$ represents the hyperbolic distance of a point of $H_2^+$ from the origin $O$ which has Cartesian coordinates $(0,1)$. We recall that $\eta$ is evaluated by means of (\ref{metric}) on the arc of a circumference with center located on the $x$ axis and joining $(x,y)$ with the origin $O$. The upper half-circumferences centered on the $x$ axis represent the geodesic lines of the space $H_2^+$ and play the same role of the straight lines in the Euclidean plane.

The angle $\alpha$ represents the slope of the tangent in $O$ to the half-circumference passing through $(x,y)$ (see Figure \ref{fig1}(a)).

\begin{figure} \label{fig1}
 \centering
 \subfigure[]
   {\includegraphics[width=5.1cm, height=6.46cm]{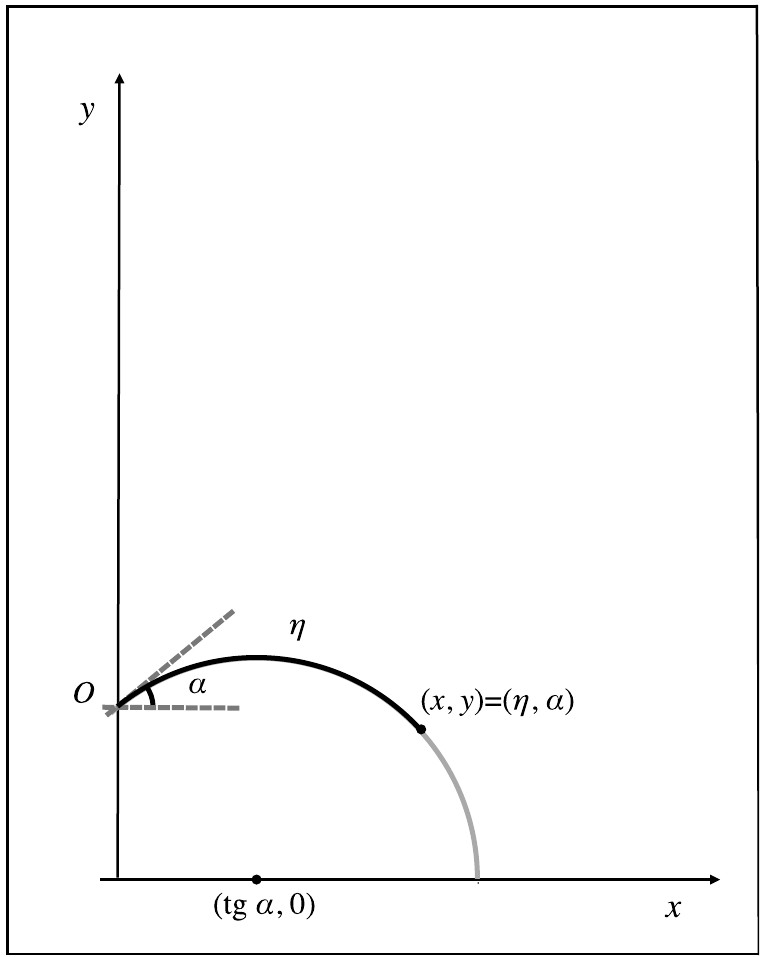}}
 \hspace{5mm}
 \subfigure[]
   {\includegraphics[width=5.95cm, height=6.46cm]{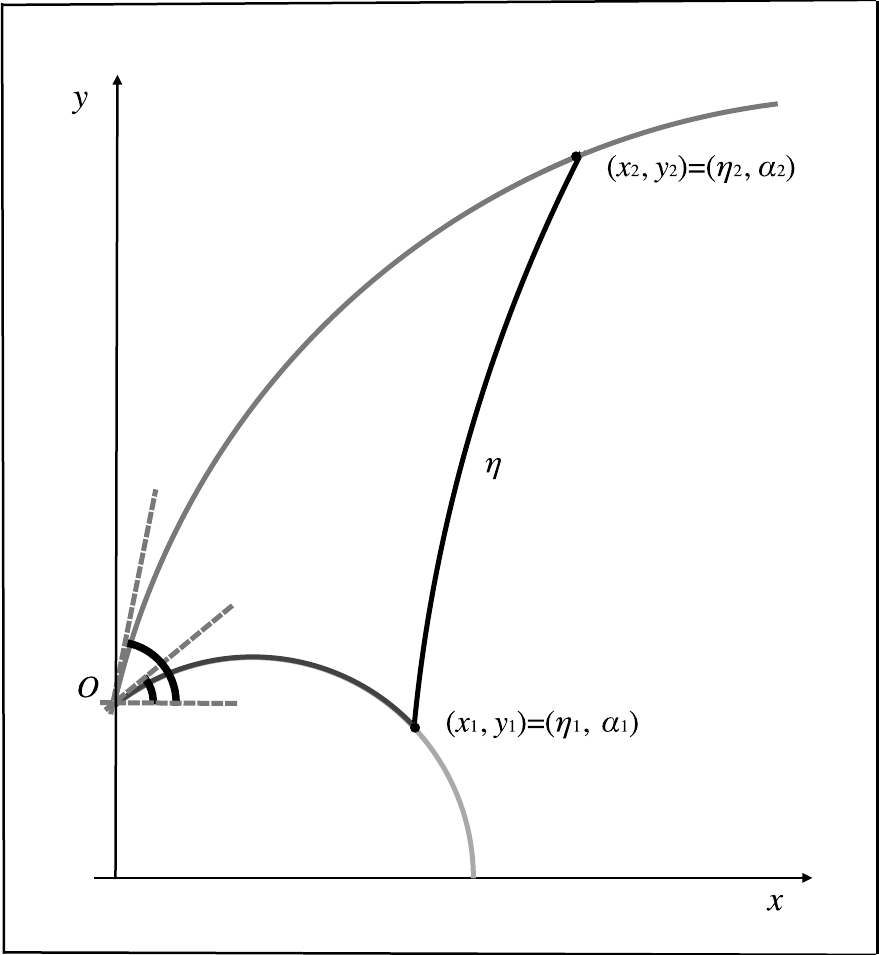}}
 \caption{Figure 1(a) illustrates the hyperbolic coordinates. Figure 1(b) refers to the hyperbolic triangle of Carnot formula. }
 \end{figure}

The formulas which relate the polar hyperbolic coordinates $(\eta,\alpha)$ to the Cartesian coordinates $(x,y)$ are (see Rogers and Williams  \cite{rog}, page 213)
\begin{equation}\label{rel}
\left\{
\begin{array}{lr} x=\frac{\sinh \eta \cos \alpha}{\cosh \eta - \sinh \eta \sin \alpha} &  \eta>0,\\
y=\frac{1}{\cosh \eta - \sinh \eta \sin \alpha} &  -\frac{\pi}{2}<\alpha<\frac{\pi}{2}.
\end{array}
\right.
\end{equation}

For each value of $\alpha$ the relevant geodesic curve is represented by the half-circumference with equation
\begin{equation}
\label{cerchio}
(x-\tan \alpha)^2+y^2=\frac{1}{\cos^2 \alpha}.
\end{equation}
For $\alpha=\frac{\pi}{2}$ we get from (\ref{cerchio}) the positive $y$ axis which also is a geodesic curve of $H_2^+$.

From (\ref{rel}) it is easy to obtain the following expression of the hyperbolic distance $\eta$ of $(x,y)$ from the origin $O$:
\begin{equation}
\label{dist}
\cosh \eta=\frac{x^2+y^2+1}{2y}.
\end{equation}
From (\ref{dist}) it can be seen that all the points having hyperbolic distance $\eta$ from the origin $O$ form a Euclidean circumference with center at $(0,\cosh \eta)$ and radius $\sinh \eta$.

The expression of the hyperbolic distance between two arbitrary points $(x_1,y_1)$ and $(x_2,y_2)$ is instead given by
\begin{equation}
\label{dist2}
\cosh \eta=\frac{(x_1-x_2)^2+y_1^2+y_2^2}{2y_1y_2}.
\end{equation}
In fact, by considering the hyperbolic triangle with vertices at $(0,1)$, $(x_1,y_1)$ and $(x_2,y_2)$, by means of the Carnot hyperbolic formula it is simple to show that the distance $\eta$ between $(x_1,y_1)$ and $(x_2,y_2)$ is given by
\begin{equation} \label{carnot}
\cosh \eta=\cosh \eta_1 \cosh \eta_2 - \sinh \eta_1 \sinh \eta_2 \cos (\alpha_1 -\alpha_2)
\end{equation}
where $(\eta_1, \alpha_1)$ and $(\eta_2, \alpha_2)$ are the hyperbolic coordinates of $(x_1,y_1)$ and  $(x_2,y_2)$, respectively (see Figure \ref{fig1}(b)). From (\ref{cerchio}) we obtain that
\begin{equation} \label{sette}
\tan \alpha_i=\frac{x_i^2+y_i^2-1}{2x_i}   \mathrm{\;\;\;\;\;\;for\;} i=1,2,
\end{equation}
and view of (\ref{dist}) and (\ref{sette}), after some calculations, formula (\ref{dist2}) appears.
 Instead of the elementary arguments of the proof above we can also invoke the group theory which reduces $(x_1, y_1)$ to $(0,1)$.

If $\alpha_1-\alpha_2=\frac{\pi}{2}$, the hyperbolic Carnot formula (\ref{carnot}) reduces to the hyperbolic Pythagorean theorem
\begin{equation}
\cosh \eta = \cosh \eta_1 \cosh \eta_2
\end{equation}
which plays an important role in the present paper.

The motion considered here is the non-Euclidean counterpart of the planar motion with orthogonal deviations studied in Orsingher \cite{O}. The main object of the investigation is the hyperbolic distance of the moving point from the origin. We are able to give explicit expressions for its mean value, also under the condition that the number of changes of direction is known. In the case of motion in $H_2^+$ with independent components (Orsingher and De Gregorio \cite{DO}) an explicit expression for the distribution of the hyperbolic distance $\eta$ has been obtained. Here, however, the components of motion are dependent and this excludes any possibility of finding the distribution of the hyperbolic distance $\eta(t)$.

We obtain the following explicit formula for the mean value of the hyperbolic distance which reads
\begin{eqnarray}
E\{ \cosh \eta(t)\}&=& e^{-\frac{\lambda t}{2}} \left\{ \cosh \frac{t}{2}\sqrt{\lambda^2+4c^2}+ \frac{\lambda}{\sqrt{\lambda^2+4c^2}} \sinh \frac{t}{2}\sqrt{\lambda^2+4c^2}   \right\} \nonumber \\
&=& Ee^{T(t)}
\end{eqnarray}
where $T(t)$ is a telegraph process with parameters $\frac{\lambda}{2}$ and $2c$.

The telegraph process represents the random motion of a particle moving with constant velocity and changing direction at Poisson-paced times (see, for example \cite{DO}).

Section $5$ is devoted to motions on the Poincar\'e half-plane where the return to the starting point is admitted and occurs at the instants of changes of direction. The mean distance from the origin of these jumping-back motions is obtained explicitly by exploiting their relationship with the motion without jumps. In the case where the return to the starting point occurs at the first Poisson event $T_1$, the mean value of the hyperbolic distance $\eta_1(t)$ reads
\begin{equation}
E\{ \cosh \eta_1(t)| N(t) \ge 1\}=\frac{\lambda}{\sqrt{\lambda^2+4c^2}} \frac{\sinh \frac{t}{2}\sqrt{\lambda^2+4c^2}}{\sinh \frac{\lambda t}{2}}.
\end{equation}

The last section considers the motion at finite velocity, with
orthogonal deviations at Poisson times, on the unit-radius sphere.
The main results concern the mean value $E\{\cos \mathrm{d}(P_0
P_t) \}$, where $\mathrm{d}(P_0 P_t)$ is the distance of the
current point $P_t$ from the starting position $P_0$. We take
profit of the analogy of the spherical motion with its counterpart
on the Poincar\'e half-plane to discuss the different situations
due to the finiteness of the space where the random motion
develops.

\section{Description of the Planar Random Motion on the Poincar\'e Half-Plane $H_2^+$}

\hspace{5mm} We start our analysis by considering a particle located at the origin $O$ of $H_2^+$. The particle initially moves on the half-circumference with center at $(0,0)$ and radius $1$. The motion of the particle develops on the geodesic lines represented by half-circles with the center located on the $x$ axis. Changes of direction are governed by a homogeneous Poisson process of rate $\lambda$.

At the occurrence of the first Poisson event, the particle starts moving on the circumference orthogonal to the previous one.

After having reached the point $P_2$, where the second Poisson event happened, the particle continues its motion on the circumference orthogonal to that joining $O$ with $P_2$ (see Figure \ref{fig:1}).

In general, at the $n$-th Poisson event, the particle is located at the point $P_n$ and starts moving on the circumference orthogonal to the geodesic curve passing through $P_n$ and the origin $O$ (consult again Figure \ref{fig:1}).

At each Poisson event the particle moves from the reached position $P$ clockwise or counter-clockwise (with probability $\frac{1}{2}$) on the circumference orthogonal to the geodesic line passing through $P$ and $O$.
\begin{figure}[t]
\centering
 \includegraphics[width=11.96cm, height=6.992cm, clip ]{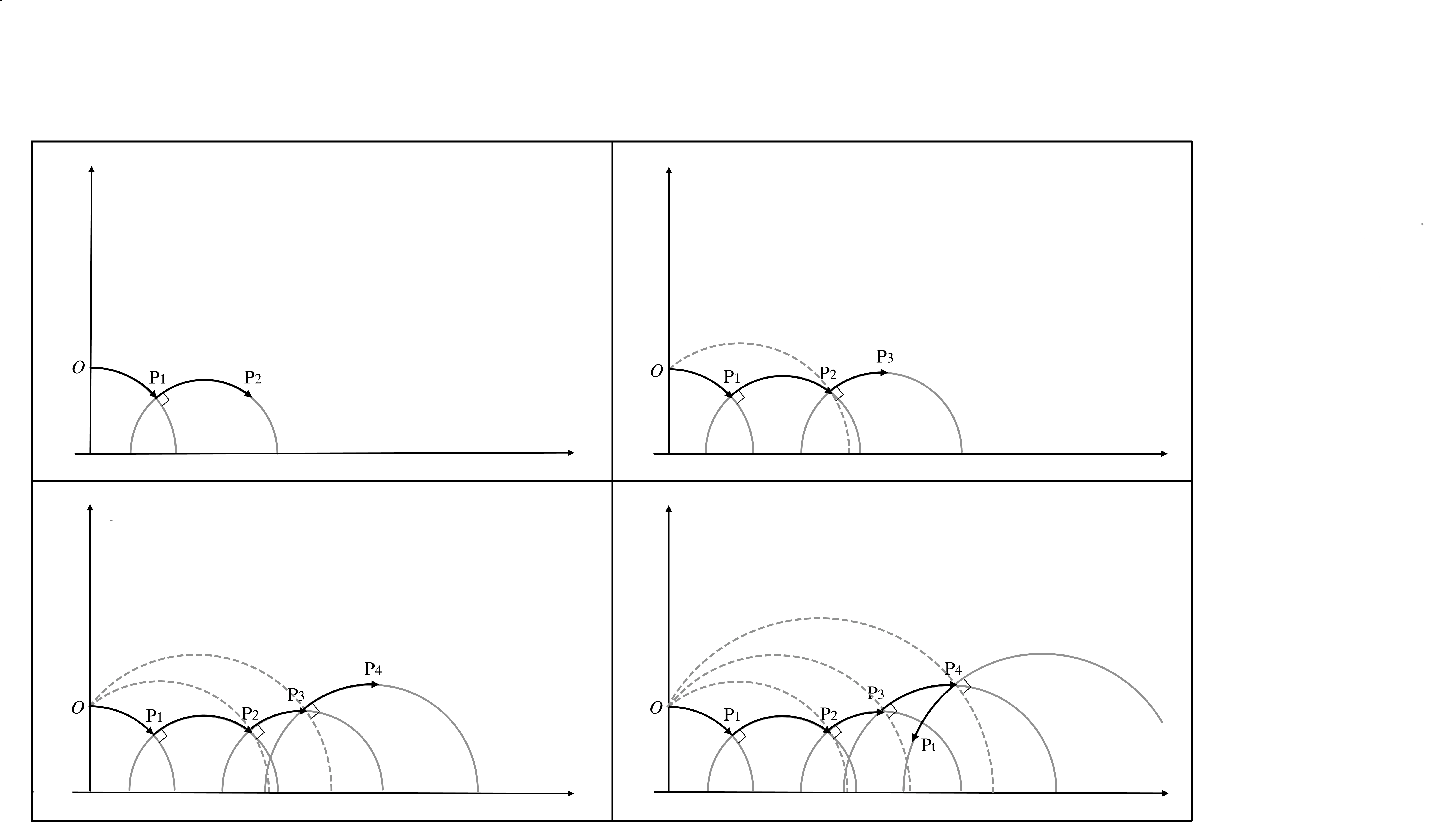}
\caption{In the first three figures a sample path where the particle chooses the outward direction is depicted. In the last one a trajectory with one step moving towards the origin is depicted.}
\label{fig:1}
\end{figure}

The hyperbolic length of the arc run by the particle during the inter-time between two successive changes of direction, occurring at $t_{k-1}$ and $t_k$ respectively, is given by $c(t_k-t_{k-1})$, with $k \ge 1$ and $t_0=0$. The velocity $c$ is assumed to be the constant hyperbolic velocity
\begin{equation}
c=\frac{\mathrm{d}s}{\mathrm{d}t}=\frac{1}{y}\sqrt{\frac{\mathrm{d}x^2+\mathrm{d}y^2}{\mathrm{d}t^2}}.
\end{equation}

The Cartesian coordinates of the points $P_k$, where the changes of direction occur, can be explicitly evaluated, but they are not important in our analysis because we study only the evolution of the hyperbolic distance from the origin of the moving particle.

The construction outlined above shows that the arcs $OP_{k-1}$, $P_{k-1}P_{k}$, and $OP_k$ form right triangles with the vertex of the right angle at $P_{k-1}$.\\
In force of the hyperbolic Pythagorean theorem we have that
\begin{equation}
\cosh \mathrm{d}(OP_{k})=\cosh \mathrm{d}(OP_{k-1}) \cosh \mathrm{d}(P_{k-1}P_{k}).
\end{equation}
The hyperbolic distance $\eta(t)$ of the moving point $P_t$ after $n$ changes of direction is thus given by
\begin{eqnarray}
\label{prod}
\cosh \eta(t)&=&\cosh \mathrm{d}(OP_t) \nonumber\\
             &=&\cosh \mathrm{d}(P_n P_t) \cosh \mathrm{d}(OP_n) \nonumber\\
             &=&\cosh c(t-t_n)\prod_{k=1}^{n}\cosh c(t_k-t_{k-1})\nonumber\\
             &=& \prod_{k=1}^{n+1} \cosh c(t_k-t_{k-1}),
\end{eqnarray}
where $t_0=0$ and $t_{n+1}=t$. The instants $t_k$, $k=0,1,\dots,n$ are uniformly distributed in the set
\begin{equation}
T=\{0=t_0<t_1< \cdots <t_k< \cdots <t_n<t_{n+1}=t\}.
\end{equation}
This means that $\cosh \eta(t)$, defined in (\ref{prod}), can be viewed as the hyperbolic distance from $O$ of the moving particle for fixed time points of the underlying Poisson process and for a fixed number $N(t)=n$ of changes of direction.

We remark that the geodesic distance (\ref{prod}) depends on how much time the particle spends on each geodesic curve (but not on the chosen direction). Of course, (\ref{prod}) depends on the number $n$ of changes of direction and on the speed $c$ of the moving particle, as well.
\begin{figure}[h]
 \centering
 \subfigure[N=2, t=6]
  {\includegraphics[width=5cm]{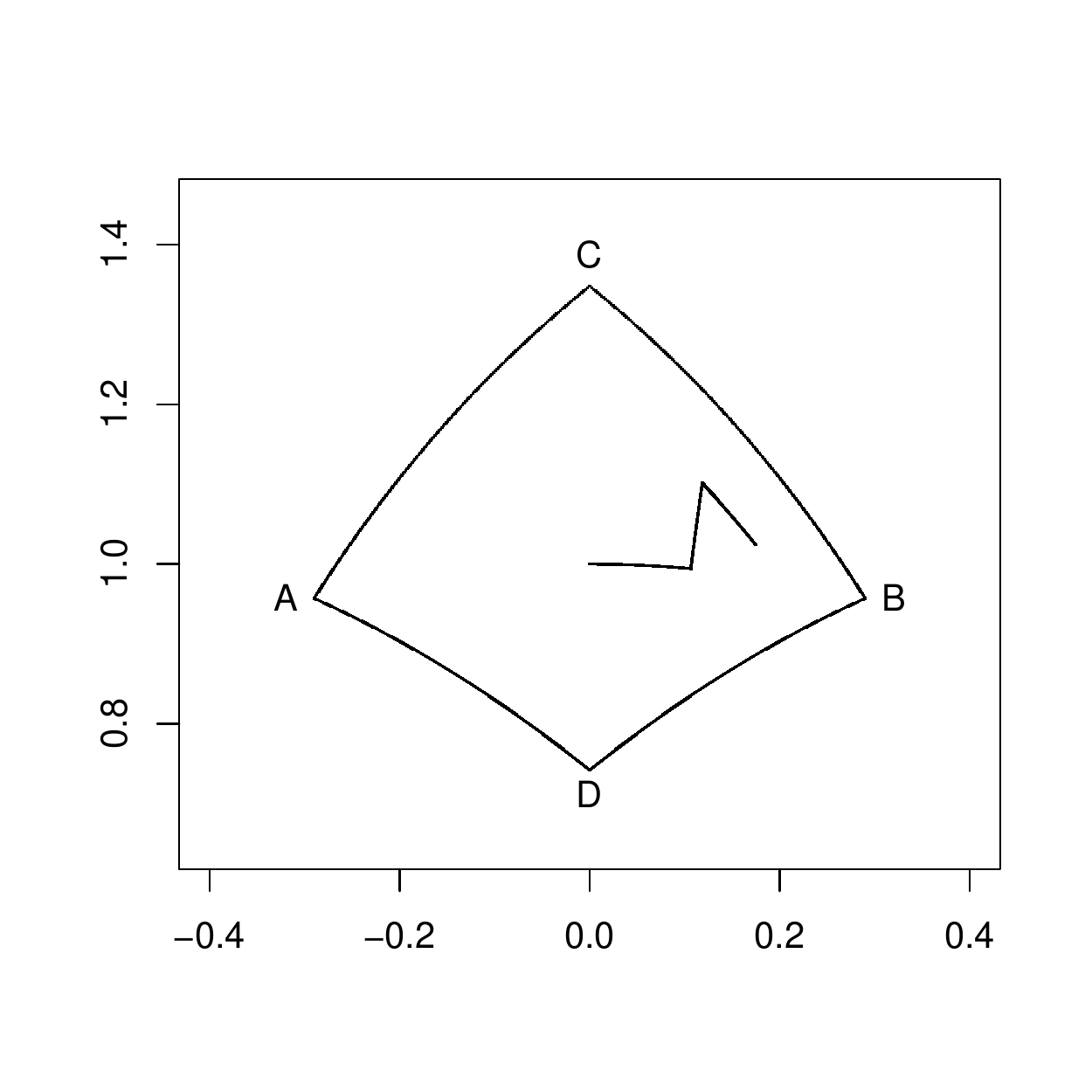}}
 \subfigure[N=2, t=25]
  {\includegraphics[width=5cm]{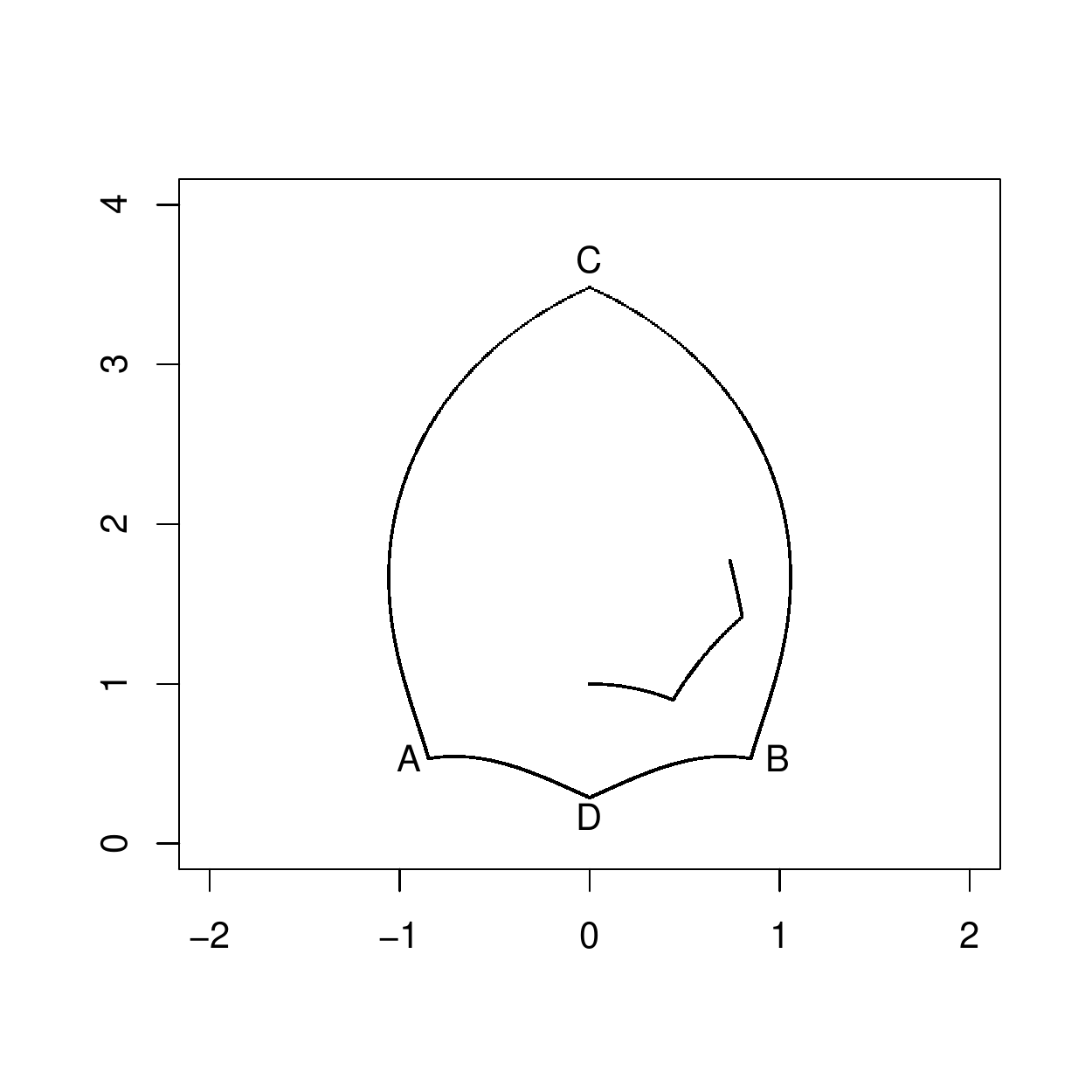}}\\ 
  \subfigure[N=2, t=40]
  {\includegraphics[width=5cm]{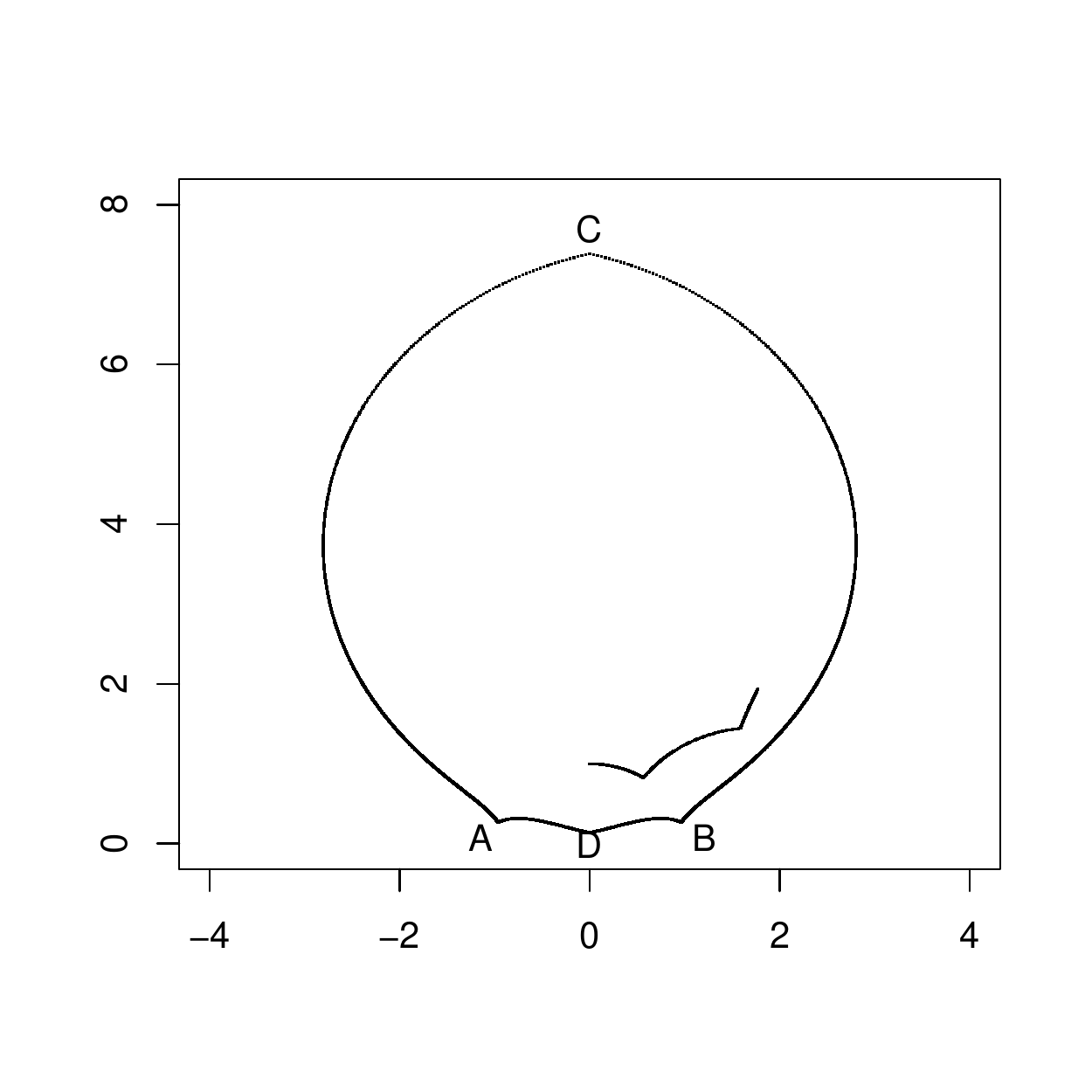}}
 \subfigure[N=2, t=50]
  {\includegraphics[width=5cm]{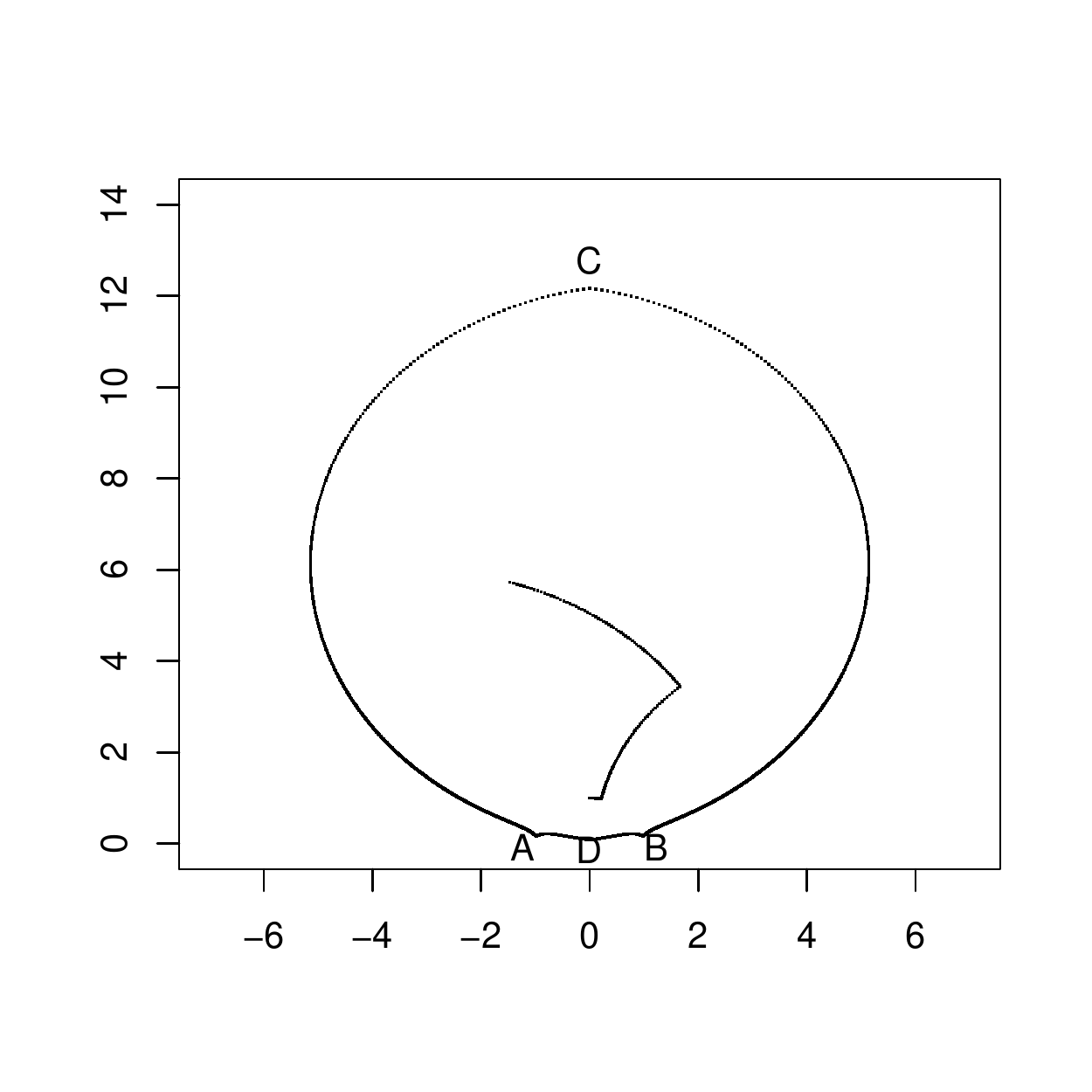}}
 \caption{The set of all possible points reachable by the process for different values of $t$  is drawn. In each domain a trajectory of the process, with $c=0.05$ and $N(t)=2$,  is depicted.}
 \label{kkk33}
 \end{figure}

The set of possible positions at different times $t$ is depicted in Figure \ref{kkk33}. The vertices $A$ and $B$ are reached when the particle never changes direction, whereas $C$ and $D$ are reached if the deviation occurs immediately after the start.

The ensemble of points having the same hyperbolic distance from $O$ at time $t$, forms the circle with center $C=(0, \cosh \eta(t))$ and radius $\sinh \eta(t)$. Since
 \begin{equation}
\cosh \eta(t)=\prod_{k=1}^{n+1}\cosh c(t_k-t_{k-1}),
\end{equation}
the ordinate of the center $C$ is obtained by successively multiplying the ordinates of the centers of equally distant points at each step. For the radius, however, such a fine interpretation is not possible (the radii do not exhibit the same multiplicative behavior) but nevertheless we will study their product
\begin{equation}  \label{ggg}
\prod_{k=1}^{n+1}\sinh c(t_k-t_{k-1})
\end{equation}
because
\begin{equation}
\sinh \eta(t) \ge  \prod_{k=1}^{n+1}\sinh c(t_k-t_{k-1})
\end{equation}
and (\ref{ggg}) represents a lower bound for the circle of equally distant points at time $t$.

\section{The Equations Related to the Mean Hyperbolic Distance}

In this section we study the conditional and unconditional mean values of the hyperbolic distance $\eta(t)$. Our first result concerns the derivation of the equations satisfied by the mean values
\begin{eqnarray} \label{bobo}
 E_n(t) &=&E\{\cosh \eta(t)| N(t)=n\} \\
 &=&\frac{n!}{t^n}\int_{0}^{t}\mathrm{d}t_1 \int_{t_1}^{t} \mathrm{d}t_2 \cdots \int_{t_{n-1}}^{t} \mathrm{d}t_n \prod_{k=1}^{n+1} \cosh c(t_k-t_{k-1})\nonumber \\
&=&\frac{n!}{t^n} I_n(t), \nonumber
\end{eqnarray}
where
\begin{equation}I_n(t)=\int_{0}^{t}\mathrm{d}t_1 \cdots \int_{t_{n-1}}^{t} \mathrm{d}t_n \prod_{k=1}^{n+1} \cosh c(t_k-t_{k-1}),\end{equation} and by
\begin{eqnarray}\label{ser}
E(t)&=&E\{\cosh \eta(t)\} \\
&=&\sum_{n=0}^{\infty} E\{\cosh\eta(t)|N(t)=n\}Pr\{N(t)=n\} \nonumber \\
&=&e^{-\lambda t}\sum_{n=0}^{\infty}\lambda^n I_n(t). \nonumber
\end{eqnarray}
At first, we state the following result concerning the evaluation of the integrals $I_n(t)$, $n \ge 1$.

\begin{lem} \label{integ}
The functions
\begin{equation} \label{integ22}
I_n(t)=\int_{0}^{t}\mathrm{d}t_1 \int_{t_1}^{t} \mathrm{d}t_2 \cdots \int_{t_{n-1}}^{t} \mathrm{d}t_n \prod_{k=1}^{n+1} \cosh c(t_k-t_{k-1}),
\end{equation}
with $t_0=0$ and $t_{n+1}=t$, satisfy the difference-differential equations
\begin{equation} \label{rrrro}
\frac{\mathrm{d}^2}{\mathrm{d}t^2}I_n=\frac{\mathrm{d}}{\mathrm{d}t}I_{n-1}+ c^2 I_n
\end{equation}
where $t >0$, $n \ge 1$, and $I_{0}(t)=\cosh{ct}$.
\end{lem}
\Dim
We first note that
\begin{eqnarray}
\frac{\mathrm{d}}{\mathrm{d}t}I_n&=&\int_{0}^{t}\mathrm{d}t_1 \cdots\int_{t_{n-2}}^{t}\mathrm{d}t_{n-1}\prod_{k=1}^{n}\cosh c(t_k-t_{k-1}) \\
&+&c\int_{0}^{t}\mathrm{d}t_1 \cdots\int_{t_{n-1}}^{t}\mathrm{d}t_n\prod_{k=1}^{n}\cosh c(t_k-t_{k-1}) \sinh c(t-t_n)\nonumber \\
&=&I_{n-1}+c\int_{0}^{t}\mathrm{d}t_1 \cdots\int_{t_{n-1}}^{t}\mathrm{d}t_n\prod_{k=1}^{n}\cosh c(t_k-t_{k-1}) \sinh c(t-t_n) \nonumber
\end{eqnarray}
and therefore
\begin{eqnarray} \label{18:53kinder}
\frac{\mathrm{d}^2}{\mathrm{d}t^2}I_n &=&\frac{\mathrm{d}}{\mathrm{d}t}I_{n-1}+c^2\int_{0}^{t}\mathrm{d}t_1 \cdots\int_{t_{n-1}}^{t}\mathrm{d}t_n\prod_{k=1}^{n+1}\cosh c(t_k-t_{k-1})\nonumber \\
&=&\frac{\mathrm{d}}{\mathrm{d}t}I_{n-1}+c^2I_n.
\end{eqnarray}
\Fine

In view of Lemma \ref{integ} we can prove also the following:

\begin{teo} \label{teo32}
The mean value $E(t)=E\{\cosh \eta(t)\}$ satisfies the second-order linear homogeneous differential equation
\begin{equation}
\frac{\mathrm{d}^2}{\mathrm{d}t^2}E(t)=-\lambda \frac{\mathrm{d}}{\mathrm{d}t}E(t)+c^2E(t)
\end{equation}
with initial conditions
\begin{equation} \label{arke}
\left\{ \begin{array}{lr}  E(0)=1, \\
\left. \frac{\mathrm{d}}{\mathrm{d}t} E(t)\right|_{t=0}=0.
\end{array}
\right.
\end{equation}\\
The explicit value of the mean hyperbolic distance is therefore
\begin{equation} \label{en}
E(t)=e^{-\frac{\lambda t}{2}}\left\{\cosh \frac{t\sqrt{\lambda^2+4c^2}}{2}+\frac{\lambda}{\sqrt{\lambda^2+4c^2}}\sinh \frac{t\sqrt{\lambda^2+4c^2}}{2} \right\}.
\end{equation}
\end{teo}
\Dim
From (\ref{ser}), it follows that
\begin{equation}\label{eq}
\frac{\mathrm{d}}{\mathrm{d}t}E(t)=-\lambda E(t)+ e^{-\lambda t} \sum_{n=0}^{\infty}\lambda^n  \frac{\mathrm{d}}{\mathrm{d}t}I_n
\end{equation}
and thus, in view of (\ref{18:53kinder}) and by letting $I_{-1}=0$, we have that
\begin{eqnarray} \label{15:58uovo}
&&\frac{\mathrm{d}^2}{\mathrm{d}t^2}E(t) \nonumber \\
&&=-\lambda\frac{\mathrm{d}}{\mathrm{d}t}E(t)-\lambda \left( \frac{\mathrm{d}}{\mathrm{d}t} E(t)+ \lambda E(t) \right)+ e^{-\lambda t}\sum_{n=0}^{\infty}\lambda^n \left( \frac{\mathrm{d}}{\mathrm{d}t}I_{n-1}+c^2 I_n \right) \nonumber \\
&&=-2 \lambda \frac{\mathrm{d}}{\mathrm{d}t}E(t)- \lambda^2 E(t)+ c^2 E(t)+ e^{-\lambda t}\sum_{n=0}^{\infty} \lambda^n \frac{\mathrm{d}}{\mathrm{d}t}I_{n-1}  \nonumber \\
&&=-2\lambda\frac{\mathrm{d}}{\mathrm{d}t}E(t)-\lambda^2 E(t)+ c^2 E(t)+\lambda \left( \frac{\mathrm{d}}{\mathrm{d}t}E(t)+\lambda E(t) \right) \nonumber \\
&&=-\lambda \frac{\mathrm{d}}{\mathrm{d}t}E(t)+c^2 E(t).
\end{eqnarray}

While it is straightforward to see that the first condition in (\ref{arke}) is verified, the second one needs some explanations: if we write
\begin{equation} \label{ago}
\left. \frac{\mathrm{d}}{\mathrm{d}t} E(t) \right|_{t=0}=\lim_{\Delta t \downarrow 0} \frac{E(\Delta t)-1 }{\Delta t}
\end{equation}
and observe that
\begin{eqnarray} \label{conduco}
&&E(\Delta t)  \\
&&=(1-\lambda \Delta t) \cosh c \Delta t + \lambda  \int_{0}^{\Delta t}\cosh ct_1 \cosh c(\Delta t -t_1) \mathrm{d}t_1+\mathrm{o}(\Delta t) \nonumber \\
&&=(1-\lambda \Delta t) \cosh c \Delta t + \frac{\lambda \Delta t}{2} \cosh c \Delta t+\frac{\lambda}{2c}\sinh{c \Delta t} +\mathrm{o}(\Delta t), \nonumber
\end{eqnarray}
by substituting (\ref{conduco}) in (\ref{ago}), the second condition emerges. The integral in (\ref{conduco}) represents the mean value $E\{ \cosh \eta(\Delta t)| N(\Delta t)=1\}$ and is in fact evaluated by applying the Pythagorean hyperbolic theorem, as in (\ref{bobo}), for $k=1$ and $t=\Delta t$.

The general solution to equation (\ref{15:58uovo}) has the form
\begin{equation} \label{gen}
E(t)=e^{-\frac{\lambda t}{2}}\left\{\mathrm{A}e^{\frac{t}{2}\sqrt{\lambda^2+4c^2}}+\mathrm{B} e^{-\frac{t}{2}\sqrt{\lambda^2+4c^2}}\right\}.
\end{equation}
By imposing the initial conditions, the constants $A$ and $B$ can be evaluated and coincide with:
\begin{equation} \label{con}
A= \frac{\lambda+\sqrt{\lambda^2+4c^2}}{2\sqrt{\lambda^2+4c^2}}, \;\;\;\;\;\;\;\;\;\;\; B=\frac{\sqrt{\lambda^2+4c^2}-\lambda}{2\sqrt{\lambda^2+4c^2}}.
\end{equation}
From (\ref{gen}) and (\ref{con}) we obtain
\begin{equation} \label{aitia}
E(t)=\frac{e^{-\frac{\lambda t}{2}}}{2}\left\{ \frac{\lambda+\sqrt{\lambda^2+4c^2}}{\sqrt{\lambda^2+4c^2}}e^{\frac{t}{2}\sqrt{\lambda^2+4c^2}}+  \frac{\sqrt{\lambda^2+4c^2}-\lambda}{\sqrt{\lambda^2+4c^2}}e^{-\frac{t}{2}\sqrt{\lambda^2+4c^2}} \right\} \nonumber
\end{equation}
so that (\ref{en}) emerges.
\Fine

\begin{oss}
The mean value $E(t)$ tends to infinity as $t \to \infty$ so that the moving particle, in the long run, either reaches the $x$ axis or moves away towards the infinity.

Of course, if $c=0$ we have that $E(t)=1$, and for $\lambda \to \infty$ we have again that $E(t)=1$ because in both cases the particle cannot leave the starting point.

If $\lambda \to 0$ we get $E(t)=\cosh ct$ because the particle will simply move on the basic geodesic line and its hyperbolic distance grows linearly with $t$.

We note that the hyperbolic distance itself tends to infinity as $t \to \infty$ because
\begin{equation} \label{22agosto}
\lim_{t \to \infty} \cosh \eta(t)= \prod_{k=1}^{\infty} \cosh \mathrm{d}(P_k P_{k-1})= \infty
\end{equation}
and (\ref{22agosto}) is the infinite product of terms bigger than one.
\end{oss}

\begin{oss}
By taking into account the difference-differential equation (\ref{rrrro}), or directly from (\ref{bobo}), it follows that the conditional mean values $E_n(t)$ satisfy the following equation with non-constant coefficients
\begin{equation} \label{17:37inglese}
\frac{\mathrm{d}^2}{\mathrm{d}t^2}E_n+\frac{2n}{t}\frac{\mathrm{d}}{\mathrm{d}t}E_n-\frac{n}{t}\frac{\mathrm{d}}{\mathrm{d}t}E_{n-1}+\frac{n^2-n}{t^2}(E_n-E_{n-1})-c^2E_{n}=0.
\end{equation}
\end{oss}

In order to obtain the explicit value of the conditional mean value $E_n(t)$ it is convenient to perform a series expansion of $E(t)$, instead of solving the difference-differential equation (\ref{17:37inglese}). In this way we can prove the following result.

\begin{teo}
The conditional mean values $E_n(t)$, $n \ge 1$, can be expressed as
\begin{eqnarray} \label{koine}
E_n(t)&=&\sum_{r=0}^{\left[\frac{n}{2} \right] }\frac{1}{2^{n}} \frac{n!}{(n-2r)!}  \sum_{j=0}^{\infty} \binom{r+j}{j} \frac{(ct)^{2j}}{(2r+2j)!} \\
&+&  \sum_{r=0}^{\left[ \frac{n-1}{2} \right]} \frac{1}{2^{n}} \frac{n!}{(n-2r-1)!}  \sum_{j=0}^{\infty} \binom{r+j}{j} \frac{(ct)^{2j}}{(2r+2j+1)!} \nonumber
\end{eqnarray}
\end{teo}
\Dim
By expanding the hyperbolic functions in (\ref{en}) we have that
\begin{eqnarray} \label{18:49naso}
E(t)&=&e^{-\frac{\lambda t}{2}} \left[  \sum_{k=0}^{\infty} \frac{1}{(2k)!} \left( \frac{t}{2}\sqrt{\lambda^2+4c^2} \right)^{2k} \right.  \\ &+& \left. \frac{\lambda}{\sqrt{\lambda^2+4c^2}} \sum_{k=0}^{\infty} \frac{1}{(2k+1)!} \left(\frac{t}{2} \sqrt{\lambda^2+4c^2} \right)^{2k+1}    \right]  \nonumber.
\end{eqnarray}
By applying the Newton binomial formula to the terms in the round brackets and by expanding $e^{\frac{\lambda t}{2}}$ it follows that
\begin{eqnarray}
E(t)& = & e^{-\lambda t} \sum_{m=0}^{\infty}\frac{1}{m!} \left( \frac{\lambda t}{2} \right)^{m}   \left[ \sum_{k=0}^{\infty} \frac{1}{(2k)!} \left(\frac{t}{2} \right)^{2k}  \sum_{r=0}^{k} \binom{k}{r} \lambda^{2r} (2c)^{2k-2r}  \right.  \nonumber \\
                  & + & \left. \lambda \sum_{k=0}^{\infty} \frac{1}{(2k+1)!} \left(\frac{t}{2} \right)^{2k+1} \sum_{r=0}^{k} \binom{k}{r} \lambda^{2r} (2c)^{2k-2r}  \right] . 
\end{eqnarray}
Finally, interchanging the summation order, it results that
\begin{eqnarray} \label{sab 10:54}
            && E(t)  \\  
                  && =  e^{-\lambda t} \left[ \sum_{m=0}^{\infty} \sum_{r=0}^{\infty} \frac{1}{m! r!} \left( \frac{\lambda t}{2} \right)^{2r+m}  \frac{(2r+m)!}{(2r+m)!}\sum_{j=0}^{\infty} \frac{(r+j)!}{j!} \frac{(ct)^{2j}}{(2r+2j)!} \right.  \nonumber \\
                  && +  \left. \sum_{m=0}^{\infty} \sum_{r=0}^{\infty} \frac{1}{m! r!} \left( \frac{\lambda t}{2} \right)^{2r+m+1}  \frac{(2r+m+1)!}{(2r+m+1)!}\sum_{j=0}^{\infty} \frac{(r+j)!}{j!} \frac{(ct)^{2j}}{(2r+2j+1)!} \right]. \nonumber
\end{eqnarray}\\
Since
\begin{equation} \label{18:43violetta}
E(t)=e^{-\lambda t} \sum_{n=0}^{\infty} \frac{(\lambda t)^n}{n!} E_n(t),
\end{equation}
from (\ref{sab 10:54}) and (\ref{18:43violetta}), we have that
\begin{eqnarray} \label{mmm}
&&E_n(t) \\
&&= \sum_{m,\,r:\;2r+m=n} \frac{1}{2^{2r+m}} \frac{(2r+m)!}{m!}  \sum_{j=0}^{\infty} \binom{r+j}{j} \frac{(ct)^{2j}}{(2r+2j)!}\nonumber \\
&&+  \sum_{m,\,r:\; 2r+m+1=n } \frac{1}{2^{2r+m+1}} \frac{(2r+m+1)!}{m!}  \sum_{j=0}^{\infty} \binom{r+j}{j} \frac{(ct)^{2j}}{(2r+2j+1)!} \nonumber \\
&&= \sum_{r=0}^{\left[\frac{n}{2} \right] }\frac{1}{2^{n}} \frac{n!}{(n-2r)!}  \sum_{j=0}^{\infty} \binom{r+j}{j} \frac{(ct)^{2j}}{(2r+2j)!} \nonumber \\
&&+  \sum_{r=0}^{\left[ \frac{n-1}{2} \right]} \frac{1}{2^{n}} \frac{n!}{(n-2r-1)!}  \sum_{j=0}^{\infty} \binom{r+j}{j} \frac{(ct)^{2j}}{(2r+2j+1)!}, \nonumber
\end{eqnarray}
and this represents the explicit form of the conditional mean values.
\Fine

\begin{oss}
We check formula (\ref{koine}) by evaluating the mean value $E_n(t)$ for $n=0,1,2,3$.

It can be noted that for $n=0$ only the term $r=0$ of the first sum in (\ref{koine}) must be considered, so that
\begin{equation}
E\{\cosh \eta(t)| N(t)=0\}=\sum_{j=0}^{\infty}\frac{(ct)^{2j}}{(2j)!}=\cosh ct.
\end{equation}

For $n=1$ both sums of (\ref{koine}) contribute to the mean value with the $r=0$ term
\begin{eqnarray}
E\{\cosh \eta(t)| N(t)=1\}&=&\frac{1}{2}\sum_{j=0}^{\infty}\frac{(ct)^{2j}}{(2j)!}+\frac{1}{2}\sum_{j=0}^{\infty}\frac{(ct)^{2j}}{(2j+1)!}\nonumber \\ &=&\frac{1}{2}\cosh ct+ \frac{1}{2ct}\sinh ct.
\end{eqnarray}

For $n=2$ we have two terms in the first sum (corresponding to $r=0,1$) and the term $r=0$ in the second sum, so that
\begin{eqnarray}
E\{\cosh \eta(t)|N(t)=2\}&=&
\frac{1}{2^2}\sum_{j=0}^{\infty}\frac{(ct)^{2j}}{(2j)!}+\frac{1}{2}\sum_{j=0}^{\infty} \binom{j+1}{j}\frac{(ct)^{2j}}{(2j+2)!}\nonumber \\
&+&\frac{1}{2}\sum_{j=0}^{\infty}\frac{(ct)^{2j}}{(2j+1)!} \nonumber \\
&=&\frac{1}{2^2}\cosh ct+ \left( \frac{1}{2^2 ct}+\frac{1}{2ct}\right)\sinh ct.
\end{eqnarray}

For $n=3$ we need to consider two terms in both sums
\begin{eqnarray}
&&E\{\cosh \eta(t)|N(t)=3\} \\
&&=\frac{1}{2^3}\sum_{j=0}^{\infty}\frac{(ct)^{2j}}{(2j)!}+\frac{3!}{2^3}\sum_{j=0}^{\infty} \binom{j+1}{j}\frac{(ct)^{2j}}{(2j+2)!} \nonumber  \\
&&+\frac{3}{2^3}\sum_{j=0}^{\infty} \frac{(ct)^{2j}}{(2j+1)!}+\frac{3}{2^2}\sum_{j=0}^{\infty}\binom{j+1}{j}\frac{(ct)^{2j}}{(2j+3)!} \nonumber \\
&&=\left(\frac{1}{2^3}+\frac{3}{2^3(ct)^2}\right)\cosh ct+ \left( \frac{6}{2^3 ct}-\frac{3}{(2ct)^3}\right)\sinh ct. \nonumber
\end{eqnarray}

The same results can be obtained directly from (\ref{bobo}) by successive integrations.
\end{oss}
For each step the ensemble of points with hyperbolic distance equal to $c(t_k-t_{k-1})$ forms a Euclidean circumference $C_k$ with radius $\sinh c(t_k-t_{k-1})$ and center located at $(0, \cosh c(t_k-t_{k-1}))$.
At time $t$, if $n$ steps have occurred, the set of points $C_t$ with hyperbolic distance equal to $\eta(t)$ is a circumference with center at $(0, \cosh \eta(t))$ and radius $\sinh \eta(t)$. Clearly \begin{equation}   \cosh \eta(t)= \prod_{k=1}^{n+1} \cosh c(t_k-t_{k-1})\end{equation} so that the ordinate of the center of $C_t$ is equal to the product of the ordinates of $C_k$. However
\begin{eqnarray}
\sinh \eta(t)&=& \sqrt{1+ \cosh^2 \eta(t)} \\
&=&\sqrt{1+\prod_{k=1}^{n+1} \cosh^2 c(t_k-t_{k-1})} \nonumber \\
&\ge& \prod_{k=1}^{n+1} \sinh c(t_k-t_{k-1}) \nonumber
\end{eqnarray}
and this shows that the quantity $\prod_{k=1}^{n+1} \sinh c(t_k-t_{k-1})$ represents a lower bound of the radius of the circle $C_t$.

\begin{teo} \label{17:02nervi}
The functions
\begin{equation}
J_n(t)=\int_{0}^{t}\mathrm{d}t_1\int_{t_1}^{t}\mathrm{d}t_2 \cdots\int_{t_{n-1}}^{t}\mathrm{d}t_n\prod_{k=1}^{n+1}\sinh c(t_k-t_{k-1}),
\end{equation}
where $t_0=0$, $t_{n+1}=t>0$, and $n\ge 1$, take the form
\begin{equation} \label{19:39orologio}
J_n(t)= \frac{t^{2n+1}c^{n+1}}{n!} \sum_{r=0}^{\infty} \frac{(n+r)!}{r!} \frac{(ct)^{2r}}{(2r+2n+1)!},
\end{equation}
where $J_0(t)=\sinh{ct}$.
\end{teo}
\Dim
We first note that the functions $J_n(t)$, $n \ge 1$, $t>0$ satisfy the difference-differential equations
\begin{equation}\label{hoho}
\frac{\mathrm{d}^2}{\mathrm{d}t^2}J_n=cJ_{n-1}+c^2J_n,  \hspace{1.2cm} n \ge 1,  \hspace{0.16cm} t>0.
\end{equation}
 Since
 \begin{equation}
\frac{\mathrm{d}}{\mathrm{d}t}J_{n}=c \int_{0}^{t}\mathrm{d}t_1 \cdots \int_{t_{n-1}}^{t} \mathrm{d}t_n \prod_{k=1}^{n} \sinh c(t_k-t_{k-1}) \cosh c(t-t_n),
\end{equation}
we have that
\begin{eqnarray}
\frac{\mathrm{d}^2}{\mathrm{d}t^2}J_{n}&=& c \int_{0}^{t}\mathrm{d}t_1 \cdots \int_{t_{n-2}}^{t} \mathrm{d}t_{n-1} \prod_{k=1}^{n} \sinh c(t_k-t_{k-1}) \nonumber \\
                                       & + &c^2 \int_{0}^{t}\mathrm{d}t_1 \cdots \int_{t_{n-1}}^{t} \mathrm{d} t_{n} \prod_{k=1}^{n+1} \sinh c(t_k-t_{k-1}) \nonumber \\
                                       & = & c J_{n-1}+c^2 J_n.
\end{eqnarray}
From (\ref{hoho}), we have that the generating function
\begin{equation} \label{19:37nervi}
G(s,t)=\sum_{n=0}^{\infty}s^n J_n
\end{equation}
satisfies the differential equation
\begin{equation}\label{lego}
\frac{\mathrm{d}^2}{\mathrm{d}t^2}G=c(s+c)G.
\end{equation}
In fact, by (\ref{hoho}), we have
\begin{equation}
\sum_{n=0}^{\infty}s^n \frac{\mathrm{d}^2}{\mathrm{d}t^2}J_n=cs\sum_{n=0}^{\infty}s^{n-1}J_{n-1}+c^2 \sum_{n=0}^{\infty} s^n J_n
\end{equation}
and this easily yields (\ref{lego}). Considering that the general solution to (\ref{lego}) is
\begin{equation} \label{sofia}
G(s,t)=\mathrm{A}e^{t \sqrt{c(s+c)}}+\mathrm{B}e^{-t\sqrt{c(s+c)}}
\end{equation}
and that $G(s, t)$ satisfies the initial conditions
\begin{equation} \label{matr}
\left\{
\begin{array}{lr}  G(s,0)=0, \\
\left. \frac{\mathrm{d}}{\mathrm{d}t}G(s,t)\right|_{t=0}=c,
\end{array}
\right.
\end{equation}
it follows that
\begin{equation} \label{dod}
G(s,t)=\frac{\sqrt{c}}{\sqrt{s+c}} \sinh t \sqrt{c(s+c)}.
\end{equation}
By expanding the $\sinh$ function in (\ref{dod}) we obtain that
\begin{eqnarray}
G(s,t)&=& \sqrt{\frac{c}{s+c}} \sum_{k=0}^{\infty} \frac{(t \sqrt{c(s+c)})^{2k+1}}{(2k+1)!} = \sum_{k=0}^{\infty} \frac{t^{2k+1}c^{k+1}(s+c)^k}{(2k+1)!} \nonumber \\
&=& \sum_{k=0}^{\infty} \sum_{j=0}^{k}\binom{k}{j} s^j c^{k-j} \frac{t^{2k+1}c^{k+1}}{(2k+1)!}= \sum_{j=0}^{\infty} s^j \left\{ \sum_{k=j}^{\infty}\binom{k}{j} c^{k-j} \frac{t^{2k+1}c^{k+1}}{(2k+1)!} \right\} \nonumber \\
&=& \sum_{j=0}^{\infty} s^j \left\{ \frac{t^{2j+1}c^{j+1}}{j!} \sum_{r=0}^{\infty}\frac{(j+r)!}{r!} \frac{(ct)^{2r}}{(2r+2j+1)!} \right\}
\end{eqnarray}
and, in view of (\ref{19:37nervi}), formula (\ref{19:39orologio}) appears.
\Fine

\begin{oss}
We consider the quantity
\begin{eqnarray} \label{sab1 11:45}
\sum_{n=0}^{\infty} \frac{n!}{t^n} J_n(t)  Pr\{N(t)=n\}&=&e^{-\lambda t} \sum_{n=0}^{\infty} \lambda^n J_n(t)=e^{-\lambda t} G(\lambda, t ) \\
&=&e^{-\lambda t} \frac{\sqrt{c}}{\sqrt{\lambda +c}} \sinh t \sqrt{c(\lambda+c)} \nonumber
\end{eqnarray}
which represents a lower bound for mean values of the radius of the circle $C$ of points with equal hyperbolic distance from the origin at time $t$. We note that the bound (\ref{sab1 11:45}) increases if
\begin{equation}
c^2+c \lambda - \lambda^2 >0.
\end{equation}
For large values of $\lambda$ the radius of the circle $C$ tends to decrease because the particle often changes direction and hardly leaves the starting point $O$.
\end{oss}

\section{About the Higher Moments of the Hyperbolic Distance}

In this section we study the conditional and unconditional higher moments of the hyperbolic distance $\eta(t)$. Our first results concern the derivation of the equations satisfied by the second-order moments
\begin{eqnarray}
M_n(t)&=&E\{\cosh^2 \eta(t)| N(t)=n\} \\
&=& \frac{n!}{t^n}\int_{0}^{t}\mathrm{d}t_1 \int_{t_1}^{t}\mathrm{d}t_2 \cdots \int_{t_{n-1}}^{t}\mathrm{d}t_n \prod_{k=1}^{n+1} \cosh^2 c(t_k-t_{k-1}) \nonumber \\
&=&\frac{n!}{t^n}U_n(t), \nonumber
\end{eqnarray}
where
\begin{equation}
U_n(t)=\int_{0}^{t}\mathrm{d}t_1 \cdots \int_{t_{n-1}}^{t}\mathrm{d}t_n \prod_{k=1}^{n+1} \cosh^2 c(t_k-t_{k-1}),
\end{equation}
and by
\begin{eqnarray} \label{10:33}
M(t)&=&E\{\cosh^2 \eta(t)\} \\
&=& \sum_{n=0}^{\infty}E\{\cosh^2 \eta(t)|N(t)=n\} Pr\{ N(t)=n\} \nonumber \\
&=&e^{-\lambda t } \sum_{n=0}^{\infty} \lambda^n U_n(t). \nonumber
\end{eqnarray}
At first, we state the following results concerning the evaluation of the integrals $U_n(t)$, $n \ge 1$.

\begin{lem} \label{integ1}
The functions
\begin{equation}\label{filia}
U_n(t)=\int_{0}^{t}\mathrm{d}t_1 \int_{t_1}^{t} \mathrm{d}t_2 \cdots \int_{t_{n-1}}^{t} \mathrm{d}t_n \prod_{k=1}^{n+1} \cosh^2 c(t_k-t_{k-1}),
\end{equation}
with $t_0=0$ and $t_{n+1}=t$, satisfy the following third-order difference-differential equations
\begin{equation} \label{ro}
\frac{\mathrm{d}^3}{\mathrm{d}t^3}U_n=\frac{\mathrm{d}^2}{\mathrm{d}t^2}U_{n-1}+ 4 c^2 \frac{\mathrm{d}}{\mathrm{d}t} U_n -2 c^2 U_{n-1},
\end{equation}
where $t>0$, $n \ge 1$, and $U_{0}(t)=\cosh^2{ct}$.
\end{lem}
\Dim
We first note that

\begin{eqnarray}\label{agape}
&&\frac{\mathrm{d}}{\mathrm{d} t} U_n \\
&&=\int_{0}^{t}\mathrm{d}t_1 \cdots \int_{t_{n-2}}^{t} \mathrm{d}t_{n-1} \prod_{k=1}^{n} \cosh^2 c(t_k-t_{k-1})\nonumber \\
&&+2c \int_{0}^{t}\mathrm{d}t_1 \cdots \int_{t_{n-1}}^{t} \mathrm{d}t_{n} \prod_{k=1}^{n} \cosh^2 c(t_k-t_{k-1}) \cosh c(t-t_n) \sinh{c(t-t_n)}\nonumber \\
&&=U_{n-1}\nonumber \\
&&+2c \int_{0}^{t}\mathrm{d}t_1 \cdots \int_{t_{n-1}}^{t} \mathrm{d}t_{n} \prod_{k=1}^{n} \cosh^2 c(t_k-t_{k-1}) \cosh c(t-t_n) \sinh{c(t-t_n)}. \nonumber
\end{eqnarray}
A further derivation yields
\begin{eqnarray}\label{kokko}
&&\frac{\mathrm{d}^2}{\mathrm{d} t^2} U_n \\
&&=\frac{\mathrm{d}}{\mathrm{d}t}U_{n-1}+2c^2 \int_{0}^{t}\mathrm{d}t_1 \cdots \int_{t_{n-1}}^{t} \mathrm{d}t_{n} \prod_{k=1}^{n+1} \cosh^2 c(t_k-t_{k-1}) \nonumber \\
&&+2c^2 \int_{0}^{t}\mathrm{d}t_1 \cdots \int_{t_{n-1}}^{t} \mathrm{d}t_{n} \prod_{k=1}^{n} \cosh^2 c(t_k-t_{k-1}) \sinh^2 c(t-t_n)\nonumber \\
&&= \frac{\mathrm{d}}{\mathrm{d}t}U_{n-1}+2c^2U_n\nonumber \\
&&+2c^2 \int_{0}^{t}\mathrm{d}t_1 \cdots \int_{t_{n-1}}^{t} \mathrm{d}t_{n} \prod_{k=1}^{n} \cosh^2 c(t_k-t_{k-1}) \sinh^2 c(t-t_n).\nonumber
\end{eqnarray}
Since it is not possible to express the integral in (\ref{kokko}) in terms of $U_n$ and its first two derivatives, a further derivation is necessary, that, in view of (\ref{agape}), leads to the following third-order difference-differential equation
\begin{eqnarray}
\frac{\mathrm{d}^3}{\mathrm{d} t^3} U_n &=& \frac{\mathrm{d}^2}{\mathrm{d}t^2}U_{n-1}+2c^2 \frac{\mathrm{d}}{\mathrm{d}t}U_n \\
&+&2^2c^3 \int_{0}^{t}\mathrm{d}t_1 \cdots \int_{t_{n-1}}^{t} \mathrm{d}t_{n} \prod_{k=1}^{n} \cosh^2 c(t_k-t_{k-1})\nonumber \\
&\cdot& \sinh c(t-t_n) \cosh c(t-t_n)\nonumber \\
&=&\frac{\mathrm{d}^2}{\mathrm{d}t^2}U_{n-1}+ 2 c^2 \frac{\mathrm{d}}{\mathrm{d}t} U_n+2c^2 \frac{\mathrm{d}}{\mathrm{d}t}U_n -2 c^2 U_{n-1}. \nonumber
\end{eqnarray}
\Fine

In view of Lemma \ref{integ1} we can prove also the following:

\begin{teo} \label{teo42}
The function $M(t)=E\{\cosh^2 \eta(t)\}$ satisfies the third-order linear differential equation
\begin{eqnarray} \label{agathos}
\frac{\mathrm{d}^3}{\mathrm{d}t^3}M(t)&=&-2\lambda\frac{\mathrm{d}^2}{\mathrm{d}t^2}
M(t)+(4c^2-\lambda^2)\frac{\mathrm{d}}{\mathrm{d}t}M(t)+2c^2\lambda M(t),
\end{eqnarray}
with initial conditions
\begin{equation}\label{kairos}
\left\{
\begin{array}{lr} M(0)=1, \\
\left. \frac{\mathrm{d}}{\mathrm{d}t}M(t)\right|_{t=0}=0,\\
\left.\frac{\mathrm{d}^2}{\mathrm{d}t^2}M(t) \right|_{t=0}=2c^2.
\end{array}
\right.
\end{equation}
\end{teo}
\Dim
By multiplying both members of (\ref{ro}) by $\lambda^n$ and summing up we have that
\begin{eqnarray}
\frac{\mathrm{d}^3}{\mathrm{d}t^3}\sum_{n=0}^{\infty}\lambda^n U_n &=&\lambda\frac{\mathrm{d}^2}{\mathrm{d}t^2}\sum_{n=1}^{\infty}\lambda^{n-1} U_{n-1}+4c^2 \frac{\mathrm{d}}{\mathrm{d}t}\sum_{n=0}^{\infty}\lambda^n U_n \nonumber \\
&-&2c^2\lambda \sum_{n=1}^{\infty}\lambda^{n-1} U_{n-1},
\end{eqnarray}
and also
\begin{equation}
\frac{\mathrm{d}^3}{\mathrm{d}t^3}\left( e^{\lambda t} M(t) \right) =\lambda \frac{\mathrm{d}^2}{\mathrm{d}t^2}\left( e^{\lambda t} M(t) \right)+ 4c^2 \frac{\mathrm{d}}{\mathrm{d}t}\left( e^{\lambda t} M(t) \right)-2c^2 \lambda e^{\lambda t} M(t),
\end{equation}
so that, after some manipulations, equation (\ref{agathos}) appears.

While the first condition in (\ref{kairos}) is obvious, the second one can be inferred from (\ref{agape}) as follows
\begin{eqnarray}
\frac{\mathrm{d}}{\mathrm{d}t}\left(e^{\lambda t}M(t) \right)&=&\frac{\mathrm{d}}{\mathrm{d}t} \sum_{n=0}^{\infty}\lambda^n U_n= \sum_{n=0}^{\infty}\lambda^n U_{n-1} \\
&+&2c\sum_{n=0}^{\infty}\lambda^n \int_{0}^{t}\mathrm{d}t_1 \cdots \int_{t_{n-1}}^{t} \mathrm{d}t_{n} \prod_{k=1}^{n} \cosh^2 c(t_k-t_{k-1}) \nonumber \\
&\cdot& \cosh c(t-t_n) \sinh{c(t-t_n)} \nonumber
\end{eqnarray}
and also
\begin{eqnarray}
\lambda e^{\lambda t} M(t)&+&e^{\lambda t}\frac{\mathrm{d}}{\mathrm{d}t} M(t)=\lambda e^{\lambda t} M(t) \nonumber \\
&+&2c \sum_{n=0}^{\infty}\lambda^n \int_{0}^{t}\mathrm{d}t_1 \cdots \int_{t_{n-1}}^{t} \mathrm{d}t_{n} \prod_{k=1}^{n} \cosh^2 c(t_k-t_{k-1})\nonumber \\
&\cdot& \cosh c(t-t_n) \sinh{c(t-t_n)}, \nonumber
\end{eqnarray}
i.e.,
\begin{equation} \label{10:44 lun}
\left. \frac{\mathrm{d}}{\mathrm{d}t} M(t)\right|_{t=0}=0
\end{equation}
since $\left. 2c \cosh ct \sinh ct \right|_{t=0}$=0. By differentiating twice (\ref{10:33}) and by taking into account (\ref{kokko}), we have that
\begin{eqnarray}
\lambda^2 e^{\lambda t} M&+&2\lambda e^{\lambda t} \frac{\mathrm{d}}{\mathrm{d}t} M +e^{\lambda t}  \frac{\mathrm{d}^2}{\mathrm{d}t^2} M= e^{\lambda t} \left( \lambda^2 M+ \lambda  \frac{\mathrm{d}}{\mathrm{d}t} M \right)+2c^2 e^{\lambda t} M\\
&+&2c^2  \sum_{n=0}^{\infty}\lambda^n \int_{0}^{t}\mathrm{d}t_1 \cdots \int_{t_{n-1}}^{t} \mathrm{d}t_{n} \prod_{k=1}^{n} \cosh^2 c(t_k-t_{k-1}) \sinh^2{c(t-t_n)},  \nonumber
\end{eqnarray}
and therefore, by considering (\ref {10:44 lun}), we obtain the second condition of (\ref{kairos}).
\Fine

In order to solve the differential equation (\ref{agathos}) we need to first solve the related third-order algebraic equation
\begin{equation}
r^3+2 \lambda r^2- (4 c^2 - \lambda^2)r - 2 c^2 \lambda=0
\end{equation}
which can be reduced to the standard form by means of the change of variable
\begin{equation} \label{apocalupto}
s=r+\frac{2 \lambda}{3}.
\end{equation}
This leads to
\begin{equation} \label{parusia}
s^3-s\left\{ \frac{\lambda^2}{3}+4 c^2 \right\}+\frac{2 \lambda}{3} \left\{ c^2- \frac{\lambda^2}{3^2} \right\}=0
\end{equation}
to which the well-known Cardano formula can be applied. In fact, for the third-order equation
\begin{equation} \label{unica via}
s^3+p s +q=0,
\end{equation}
the solution can be expressed as
\begin{equation} \label{cardano}
s=\sqrt[3]{-\frac{q}{2}+ \sqrt{\frac{p^3}{3^3}+\frac{q^2}{2^2}}}+\sqrt[3]{-\frac{q}{2}- \sqrt{\frac{p^3}{3^3}+\frac{q^2}{2^2}}}.
\end{equation}
By comparing (\ref{parusia}) and (\ref{unica via}) it results
\begin{equation}
\frac{p^3}{3^3}+\frac{q^2}{2^2}=-\frac{c^2}{3^3}\left[ (2^3c^2+\lambda^2)^2 + \lambda^2(\lambda^2-3c^2) \right],
\end{equation}
\begin{equation}
-\frac{q}{2}=-\frac{\lambda}{3}\left( c^2- \frac{\lambda^2}{3^2}\right).
\end{equation}
The simplest case is that of $c=\frac{\lambda}{3}$ for which the solutions of  (\ref{parusia}) are $s_1=0$, $s_2=\sqrt{7}c$ and $s_3=-\sqrt{7}c$. After some calculations we get that
\begin{equation}
E\{\cosh^2 \eta(t) \}=\frac{e^{-2ct}}{7}\left\{1+6 \cosh \sqrt{7}ct+2 \sqrt{7} \sinh \sqrt{7}ct  \right\}.
\end{equation}

Following Lemma (\ref{integ1}) we can prove a more general result:

\begin{teo}
The functions
\begin{equation}
K^m_n(t)=\int_{0}^{t}\mathrm{d}t_1 \cdots \int_{t_{n-1}}^{t}\mathrm{d}t_n \prod_{k=1}^{n+1} \cosh^m c(t_k-t_{k-1}),
\end{equation}
with $t_0=0$ and $t_{n+1}=t$, are solutions of difference-differential equations of order $m+1$.
\end{teo}
\Dim
For $m=1$ and $m=2$ this statement has already been shown above since, in Theorem \ref{teo32} and Theorem \ref{teo42}, we have obtained that
\begin{equation}
\frac{\mathrm{d}^2}{\mathrm{d}t^2} K^1_{n}+\lambda \frac{\mathrm{d}}{\mathrm{d}t} K^1_{n-1}-c^2 K_n^1 =0,
\end{equation}
and
\begin{equation}
\frac{\mathrm{d}^3}{\mathrm{d}t^3} K^2_{n}-\frac{\mathrm{d}^2}{\mathrm{d}t^2} K^2_{n-1}-4 c^2 \frac{\mathrm{d}}{\mathrm{d}t} K_n^2+2 c^2 K_{n-1}^2 =0.
\end{equation}\\
We easily see that
\begin{eqnarray} \label{deloi}
\frac{\mathrm{d}}{\mathrm{d}t} K^m_{n}&=& K^m_{n-1} \\ &+& c\, m \int_{0}^{t} \mathrm{d}t_1 \cdots \int_{t_{n-1}}^{t} \mathrm{d}t_n \prod_{k=1}^{n} \cosh^m c(t_k-t_{k-1}) \nonumber \\
&\cdot&\cosh^{m-1}c(t-t_n) \sinh c(t-t_n), \nonumber
\end{eqnarray}
and
\begin{eqnarray} \label{talassa}
\frac{\mathrm{d}^2}{\mathrm{d}t^2} K^m_{n}&=& \frac{\mathrm{d}}{\mathrm{d}t} K^m_{n-1}+ c^2 m K_n^{m} \\
 &+& c^2 m (m-1) \int_{0}^{t} \mathrm{d}t_1 \cdots \int_{t_{n-1}}^{t} \mathrm{d}t_n \prod_{k=1}^{n} \cosh^m c(t_k-t_{k-1}) \nonumber \\
&\cdot& \cosh^{m-2}c(t-t_n) \sinh^2 c(t-t_n). \nonumber
\end{eqnarray}
In view of (\ref{deloi}) it also results
\begin{eqnarray} \label{gulaxsa}
\frac{\mathrm{d}^3}{\mathrm{d}t^3} K^m_{n}&=& \frac{\mathrm{d}^2}{\mathrm{d}t^2} K^m_{n-1}+ c^2 m  \frac{\mathrm{d}}{\mathrm{d}t} K_n^{m} + 2 c^2 (m-1) \left\{ \frac{\mathrm{d}}{\mathrm{d}t}K_{n}^m-K_{n-1}^m \right\} \nonumber \\
&+&c^3 m (m-1) (m-2) \int_{0}^{t} \mathrm{d}t_1 \cdots \int_{t_{n-1}}^{t} \mathrm{d}t_n \prod_{k=1}^{n} \cosh^m c(t_k-t_{k-1}) \nonumber \\
&\cdot& \cosh^{m-3}c(t-t_n) \sinh^3 c(t-t_n).
\end{eqnarray}\\
After $(m-1)$ derivatives the following equation is obtained
\begin{eqnarray} \label{deloo}
\frac{\mathrm{d}^{m-1}}{\mathrm{d}t^{m-1}} K^m_{n}&=& \frac{\mathrm{d}^{m-2}}{\mathrm{d}t^{m-2}} K^m_{n-1}+ c^2 m \frac{\mathrm{d}^{m-3}}{\mathrm{d}t^{m-3}} K_n^{m} + \cdots  + \nonumber \\
&&+ c^{m-1} m (m-1) \cdots (m-(m-1)+1) \nonumber \\
 &\cdot&    \int_{0}^{t} \mathrm{d}t_1 \cdots \int_{t_{n-1}}^{t} \mathrm{d}t_n \prod_{k=1}^{n} \cosh^m c(t_k-t_{k-1})\nonumber \\
 &\cdot & \cosh c(t-t_n) \sinh^{m-1} c(t-t_n),
\end{eqnarray}
and the next derivative gives
\begin{eqnarray} \label{manthano}
&&\frac{\mathrm{d}^{m}}{\mathrm{d}t^{m}} K^m_{n}= \frac{\mathrm{d}^{m-1}}{\mathrm{d}t^{m-1}} K^m_{n-1}+ c^2 m \frac{\mathrm{d}^{m-2}}{\mathrm{d}t^{m-2}} K_n^{m} + \cdots  +\nonumber \\
&+& c^{m} m (m-1) \cdots 2 \nonumber \\
&\cdot& \int_{0}^{t} \mathrm{d}t_1 \cdots \int_{t_{n-1}}^{t} \mathrm{d}t_n \prod_{k=1}^{n} \cosh^m c(t_k-t_{k-1}) \sinh^{m} c(t-t_n) \nonumber \\
 &+& c^{m} m (m-1) \cdots 2 \cdot (m-1) \nonumber \\
 &\cdot& \int_{0}^{t} \mathrm{d}t_1 \cdots \int_{t_{n-1}}^{t} \mathrm{d}t_n \prod_{k=1}^{n} \cosh^m c(t_k-t_{k-1})   \nonumber \\
 &\cdot& \cosh^2 c(t-t_n) \sinh^{m-2} c(t-t_n) .
\end{eqnarray}
The second integral of (\ref{manthano}) can be expressed in terms of the derivatives of order $(m-2)$ and lower. \\By further differentiating equation (\ref{manthano}) it turns out that, because of (\ref{deloo}), the derivative of the first integral in (\ref{manthano}) can be expressed in terms of the derivatives of order ($m-1$) and lower. The theorem is thus proved.
\Fine

Likewise  Theorem \ref{17:02nervi}, the following theorem holds:

\begin{teo}
The function
\begin{equation}
V_{n}(t)=\int_{0}^{t}\mathrm{d}t_1 \int_{t_1}^{t} \mathrm{d}t_2 \cdots \int_{t_{n-1}}^{t}\mathrm{d}t_n \prod_{k=1}^{n+1} \sinh^2 c(t_k -t_{k-1})
\end{equation}
with $t_0=0$ and $t_{n+1}=t$, satisfies the third-order difference-differential equation
\begin{equation}
\frac{\mathrm{d}^3}{\mathrm{d}t^3}V_{n}=4 c^2 \frac{\mathrm{d}}{\mathrm{d}t}V_{n}+2c^2 V_{n-1}
\end{equation}
where $t>0$, $n \ge 1$, and $V_0(t)=\sinh^2 ct$.
\end{teo}
\Dim
We first note that
\begin{equation}\label{dios}
\frac{\mathrm{d}}{\mathrm{d}t} V_{n}=2c \int_{0}^{t}\mathrm{d}t_1 \cdots \int_{t_{n-1}}^{t}\mathrm{d}t_n \prod_{k=1}^{n} \sinh^2 c(t_k -t_{k-1}) \sinh c(t-t_n) \cosh c(t-t_n)
\end{equation}
and therefore
\begin{equation}
\frac{\mathrm{d}^2}{\mathrm{d}t^2} V_{n}=2c^2 V_n+ 2c^2 \int_{0}^{t}\mathrm{d}t_1 \cdots \int_{t_{n-1}}^{t}\mathrm{d}t_n \prod_{k=1}^{n} \sinh^2 c(t_k -t_{k-1}) \cosh^2 c(t-t_n),
\end{equation}
and
\begin{eqnarray}\label{zeus}
\frac{\mathrm{d}^3}{\mathrm{d}t^3} V_{n}&=&2 c^2 \frac{\mathrm{d}}{\mathrm{d}t} V_n+ 2 c^2 V_{n-1}\nonumber \\ &+&4 c^3 \int_{0}^{t}\mathrm{d}t_1 \cdots \int_{t_{n-1}}^{t}\mathrm{d}t_n \prod_{k=1}^{n} \sinh^2 c(t_k -t_{k-1})\nonumber \\ &\cdot& \sinh c(t-t_n) \cosh c(t-t_n).
\end{eqnarray}
Finally, by substituting (\ref{dios}) in (\ref{zeus}), we obtain
\begin{equation}
\frac{\mathrm{d}^3}{\mathrm{d}t^3}V_{n}=4 c^2 \frac{\mathrm{d}}{\mathrm{d}t}V_{n}+2c^2 V_{n-1}.
\end{equation}
\Fine


\section{Motions with Jumps Backwards to the Starting Point}

We here examine the planar motion dealt with so far assuming now that, at the instants of changes of direction, the particle can return to the starting point  and commence its motion from scratch.

The new motion and the original one are governed by the same Poisson process so that changes of direction occur simultaneously in the original as well as in the new motion starting afresh from the origin. This implies that the arcs of the original sample path and those of the new trajectories have the same hyperbolic length. However, the angles formed by successive segments differ in order to make the hyperbolic Pythagorean theorem applicable to the trajectories of the new motion.

In order to make our description clearer, we consider the case where, in the interval $(0,t)$, $N(t)=n$ Poisson events ($n \ge 1$) occur and we assume that the jump to the origin happens at the first change of direction, i.e., at the instant $t_1$.
The instants of changes of direction for the new motion are \begin{equation} t'_k=t_{k+1}-t_1\end{equation} where $k=0, \cdots, n$ with $t'_0=0$ and $t'_n=t-t_1$ and the hyperbolic lengths of the corresponding arcs are
\begin{equation}
c(t'_{k}-t'_{k-1})=c(t_{k+1}-t_k).
\end{equation}
Therefore, at the instant $t$, the hyperbolic distance from the origin of the particle performing the motion which has jumped back to $O$ at time $t_1$ is
\begin{eqnarray} \label{gallia}
\prod_{k=1}^{n} \cosh c(t'_k-t'_{k-1}) &=& \prod_{k=1}^{n} \cosh c(t_{k+1}-t_k) \nonumber \\
&=& \prod_{k=2}^{n+1} \cosh c(t_k -t_{k-1})
\end{eqnarray}
where $0=t'_0<t'_1< \cdots < t'_n=t-t_1$ and $t_{k+1}=t'_k+t_1$. Formula (\ref{gallia}) shows that the new motion has an hyperbolic distance equal to that of the original motion where the first step has been deleted. However, the distance between the position $P_t$ and the origin $O$ of the moving particle which jumped back to $O$ after having reached the position $P_1$, is different from the distance of $P_t$ from $P_1$ since the angle between successive steps must be readjusted in order to apply the hyperbolic Pythagorean theorem.

If we denote by $T_1$ the random instant of the return to the starting point (occurring at the first Poisson event), we have that
\begin{eqnarray} \label{cisalpina}
&&E\{ \cosh \eta_1(t)I_{\{ N(t) \ge 1\}} | N(t)=n\}\\ &&=E\{ \cosh \eta (t-T_1) I_{\{T_1 \le t \}} | N(t)=n\}  \nonumber \\
&&=\int_{0}^{t} E\{ \cosh \eta (t-T_1) I_{\{ T_1 \in \mathrm{d}t_1 \}} | N(t)=n \}  \mathrm{d} t_1 \nonumber \\
&&=\int_{0}^{t}  E\{ \cosh \eta (t-T_1) |T_1=t_1, N(t)=n \}   Pr\{T_1 \in \mathrm{d} t_1 | N(t)=n\} \mathrm{d} t_1. \nonumber
\end{eqnarray}
By observing that
\begin{eqnarray}
E\{ \cosh \eta (t-T_1)|T_1=t_1, N(t)=n \}&=&E\{ \cosh \eta (t-t_1) | N(t)=n-1\} \nonumber \\
&=& \frac{(n-1)!}{(t-t_1)^{n-1}} I_{n-1}(t-t_1),
\end{eqnarray}
and that
\begin{equation}
Pr\{T_1 \in \mathrm{d} t_1 |N(t)=n\}=\frac{n!}{t^n} \frac{(t-t_1)^{n-1}}{(n-1)!} \mathrm{d}t_1
\end{equation}
with $0<t_1<t$, formula (\ref{cisalpina}) becomes
\begin{equation} \label{17:31follia}
E\{ \cosh \eta_1(t) I_{\{ N(t) \ge 1\}} |N(t)=n\}= \frac{n!}{t^n} \int_{0}^{t} I_{n-1}(t-t_1) \mathrm{d}t_1.
\end{equation}
From (\ref{17:31follia}) we have that the mean hyperbolic distance for the particle which returns to $O$ at time $T_1$ has the form:
\begin{eqnarray}
E\{\cosh \eta_1(t) | N(t) \ge 1\}&=&\frac{e^{-\lambda t}}{Pr\{N(t) \ge 1\}} \sum_{n=1}^{\infty} \lambda^n \int_{0}^{t} I_{n-1}(t-t_1)\mathrm{d}t_1 \nonumber \\
&=&\frac{\lambda  e^{-\lambda t} }{Pr\{N(t) \ge 1\}}\int_{0}^{t} e^{\lambda (t-t_1)} E(t-t_1) \mathrm{d}t_1
\end{eqnarray}

We give here a general expression for the mean value of the hyperbolic distance of a particle which returns to the origin for the last time at the $k$-th Poisson event $T_k$. We shall denote the distance by the following equivalent notation $\eta(t-T_k)=\eta_k(t)$ where the first expression underlines that the particle starts from scratch at time $T_k$ and then moves away for the remaining interval of length $t-T_k$.
In the general case we have the result stated in the next theorem:

\begin{teo}
If $N(t) \ge k$, then the mean value of the hyperbolic distance $\eta_k$ is equal to
\begin{eqnarray} \label{9:36coppa}
&&E\{ \cosh \eta_k(t) | N(t) \ge k\} \\
&&= \frac{\lambda^k e^{-\lambda t} }{Pr\{N(t) \ge k\}}  \int_{0}^{t} \mathrm{d}t_1 \cdots \int_{t_{k-1}}^{t}  e^{\lambda(t- t_k)} E(t-t_k) \mathrm{d}t_k \nonumber \\
&&=\frac{\lambda^k e^{-\lambda t}}{Pr\{N(t) \ge k\} (k-1)!} \int_{0}^{t} e^{\lambda(t- t_k)} t_k^{k-1} E(t-t_k) \mathrm{d}t_k, \nonumber
\end{eqnarray}
where $E(t)$ is given by (\ref{en}).
\end{teo}
\Dim
We start by observing that
\begin{eqnarray} \label{13:49 tacchino}
&&E\{ \cosh \eta_k(t) | N(t) \ge k\} \\
&&= \sum_{n=k}^{\infty} E\{ \cosh \eta_k(t) I_{\{N(t) =n\}} | N(t) \ge k \} \nonumber \\
&&=\sum_{n=k}^{\infty} E\{ \cosh \eta_k(t) I_{\{ N(t) \ge k\}} | N(t)=n\} \frac{Pr\{N(t)=n\}}{Pr\{N(t) \ge k\}} \nonumber \\
&&= \sum_{n=k}^{\infty} E\{ \cosh \eta_k(t) I_{\{ N(t) \ge k\}} | N(t)=n\} Pr\{N(t)=n | N(t) \ge k\}. \nonumber
\end{eqnarray}
Since $T_k= \inf\{t: N(t)=k\}$, the conditional mean value inside the sum can be developed as follows
\begin{eqnarray}
&&E\{ \cosh \eta_k(t) I_{\{N(t) \ge k\}} | N(t)=n \} \\
&&= E\{ \cosh \eta(t-T_k) I_{\{ T_k \le t\}} |N(t)=n\} \nonumber \\
&&=\int_{0}^{t}E\{\cosh \eta(t-t_k)I_{\{T_k \in \mathrm{d}t_k\} } |N(t)=n \} \mathrm{d}t_k\nonumber \\
&&=\int_{0}^{t}E\{\cosh \eta(t-t_k) | T_k =t_k, N(t)=n \} Pr\{ T_k \in \mathrm{d}t_k | N(t)=n\} \mathrm{d}t_k. \nonumber
\end{eqnarray}
In view of (\ref{bobo}) we have that
\begin{eqnarray}
E\{ \cosh \eta(t-T_k) |T_k=t_k, N(t)=n\}&=&E\{ \cosh \eta(t-t_k)| N(t-t_k)=n-k\} \nonumber \\
&=& \frac{(n-k)!}{(t-t_k)^{n-k}} I_{n-k}(t-t_k),
\end{eqnarray}
and on the base of well-known properties of the Poisson process we have that
\begin{equation}
Pr\{ T_k \in \mathrm{d} t_k | N(t)=n\} = \frac{n!}{t^n} \frac{(t-t_k)^{n-k}}{(n-k)!} \frac{t_k^{k-1}}{(k-1)!} \mathrm{d}t_k
\end{equation}
where $0<t_k<t$. In conclusion we have that
\begin{equation}
E\{ \cosh \eta_k(t) I_{\{N(t) \ge k\}} | N(t)=n\}= \frac{n!}{t^n}\frac{1}{(k-1)!} \int_{0}^{t} t_k^{k-1} I_{n-k}(t-t_k) \mathrm{d} t_k
\end{equation}
and, from this and (\ref{13:49 tacchino}), it follows that
\begin{eqnarray}
&&E\{ \cosh \eta_k(t) | N(t) \ge k\} \\
&&= \sum_{n=k}^{\infty}  \frac{n!}{t^n(k-1)!} \int_{0}^{t} t_k^{k-1} I_{n-k}(t-t_k) \mathrm{d} t_k \frac{e^{-\lambda t} (\lambda t)^n}{n! Pr\{N(t) \ge k\}} \nonumber \\
&&= \frac{ \lambda^k e^{-\lambda t} }{Pr\{N(t)\ge k\}(k-1)!} \int_{0}^{t} e^{\lambda (t-t_k)} t_k^{k-1} E(t-t_k) \mathrm{d}t_k. \nonumber
\end{eqnarray}
Finally, in view of Cauchy formula of multiple integrals, we obtain that
\begin{eqnarray}
&&\frac{ \lambda^k e^{-\lambda t} }{Pr\{N(t)\ge k\}(k-1)!} \int_{0}^{t}  e^{\lambda(t- t_k)} t_k^{k-1} E(t-t_k) \mathrm{d}t_k \\
&&=\frac{\lambda^k e^{-\lambda t} }{Pr\{N(t) \ge k\}}  \int_{0}^{t} \mathrm{d}t_1 \cdots \int_{t_{k-1}}^{t} e^{\lambda(t- t_k)} E(t-t_k) \mathrm{d}t_k. \nonumber
\end{eqnarray}
\Fine


\begin{teo}
The mean of the hyperbolic distance of the moving particle returning to the origin at the $k$-th change of direction is
\begin{eqnarray}
&&E\{ \cosh \eta_k(t) | N(t) \ge k \} \nonumber \\
&&=\frac{\lambda^k e^{-\lambda t}}{\sqrt{\lambda^2+4c^2} Pr\{N(t) \ge k\}} \left\{ \frac{e^{At}}{A^{k-1}}- \frac{e^{Bt}}{B^{k-1}}+\sum_{i=1}^{k-1} \left(  \frac{1}{B^i}-\frac{1}{A^i}\right) \frac{t^{k-i-1}}{(k-i-1)!}   \right\} \nonumber
\end{eqnarray}\\
 where
 \begin{equation}  \label{arkee} A=\frac{1}{2} (\lambda + \sqrt{\lambda^2+4c^2} ), \;\;\;\;\;\;\;\;\;\;\;\;  B=\frac{1}{2}(\lambda - \sqrt{\lambda^2+4c^2} ). \end{equation}
For $k=1$, the sum in (\ref{arkee}) is intended to be zero.
\end{teo}
\Dim
We can prove (\ref{arkee}) by applying both formulas in (\ref{9:36coppa}). We start our proof by employing the first one:
\begin{eqnarray} \label{9:47burro}
&&E\{ \cosh \eta_k(t) | N(t) \ge k \} \\
&&= \frac{\lambda^k e^{-\lambda t} }{Pr\{N(t) \ge k\}}  \int_{0}^{t} \mathrm{d}t_1 \cdots \int_{t_{k-1}}^{t}  e^{\lambda (t-t_k)} E(t-t_k) \mathrm{d}t_k. \nonumber
\end{eqnarray}
Therefore, in view of (\ref{aitia}), formula (\ref{9:47burro}) becomes
\begin{eqnarray}
&&E\{ \cosh \eta_k(t) | N(t) \ge k  \}  \nonumber \\
&&= \frac{\lambda^ke^{-\lambda t}}{Pr\{N(t) \ge k\}} \int_{0}^{t} \mathrm{d}t_1 \cdots \int_{t_{k-1}}^{t}  e^{\lambda(t-t_k)} \left\{ \frac{e^{-\frac{\lambda}{2}(t-t_k)}}{2}  \left[  \left( \frac{\lambda+\sqrt{\lambda^2+4c^2}}{\sqrt{\lambda^2+4c^2}} \right) \cdot \right. \right.  \nonumber \\
&& \left.  e^{\frac{(t-t_k)}{2} \sqrt{\lambda^2+4c^2}}   + \left. \left( \frac{-\lambda+\sqrt{\lambda^2+4c^2}}{\sqrt{\lambda^2+4c^2}} \right) e^{-\frac{(t-t_k)}{2} \sqrt{\lambda^2+4c^2}}   \right]  \right\} \mathrm{d}t_k . \nonumber
\end{eqnarray}
By introducing $A$ and $B$ as in (\ref{arkee}), we can easily determine the $k$-fold integral
\begin{eqnarray}
&&E\{ \cosh \eta_k(t) | N(t) \ge k  \}  \\
&&=\frac{\lambda^ke^{-\lambda t}}{\sqrt{\lambda^2+4c^2} Pr\{N(t) \ge k\}} \int_{0}^{t} \mathrm{d}t_1 \cdots \int_{t_{k-1}}^{t} \left\{ A e^{A(t-t_k)} - B e^{B(t-t_k)} \right\}  \mathrm{d}t_k  \nonumber \\
&&=\frac{\lambda^ke^{-\lambda t}}{\sqrt{\lambda^2+4c^2} Pr\{N(t) \ge k\}} \int_{0}^{t} \mathrm{d}t_1 \cdots \int_{t_{k-2}}^{t} \left\{ e^{A(t-t_{k-1})}-e^{B(t-t_{k-1})} \right\}  \mathrm{d}t_{k-1}  \nonumber \\
&&=\frac{\lambda^k e^{-\lambda t}}{\sqrt{\lambda^2+4c^2} Pr\{N(t) \ge k\}} \int_{0}^{t} \mathrm{d}t_1 \cdots \int_{t_{k-3}}^{t} \left\{ \frac{e^{A(t-t_{k-2})}}{A}-\frac{e^{B(t-t_{k-2})}}{B}+\frac{1}{B}-\frac{1}{A}   \right\}  \mathrm{d}t_{k-2} . \nonumber
\end{eqnarray}
At the $j$-th stage the integral becomes
\begin{eqnarray}
&&E\{ \cosh \eta_k(t) | N(t) \ge k \}  \nonumber \\
 &&=\frac{\lambda^ke^{-\lambda t}}{\sqrt{\lambda^2+4c^2} Pr\{N(t) \ge k\}} \int_{0}^{t} \mathrm{d}t_1 \cdots \int_{t_{k-j-1}}^{t} \left\{ \frac{e^{A(t-t_{k-j})}}{A^{j-1}}-\frac{e^{B(t-t_{k-j})}}{B^{j-1}} \right. \nonumber \\
 &&\hspace{5cm} \left. + \sum_{i=1}^{j-1}\left(\frac{1}{B^i}-\frac{1}{A^i}\right)\frac{(t-t_{k-j})^{j-i-1}}{(j-i-1)!}\right\} .  \nonumber
 \end{eqnarray}
 At the $k-1$-th stage the integral becomes
 \begin{eqnarray}
 &&E\{ \cosh \eta_k(t) | N(t) \ge k \}   \\
 &&=\frac{\lambda^k e^{-\lambda t}}{\sqrt{\lambda^2+4c^2} Pr\{N(t) \ge k\}} \int_{0}^{t} \mathrm{d}t_1 \left\{  \frac{e^{A(t-t_1)}}{A^{k-2}}- \frac{e^{B(t-t_1)}}{B^{k-2}}\right. \nonumber \\
 &&\hspace{4cm} \left.+\sum_{i=1}^{k-2}\left( \frac{1}{B^i}- \frac{1}{A^i}\right) \frac{ (t-t_1)^{k-i-2}}{(k-i-2)!}\right\}. \nonumber
 \end{eqnarray}
At the $k$-th integration we obtain formula (\ref{arkee}).

By means of the second formula in (\ref{9:36coppa}) and by repeated integrations by parts we can obtain again result (\ref{arkee}).
\Fine

\begin{oss}
For $k=1$ we have that
\begin{equation} \label{transalpina}
E\{ \cosh \eta_1(t) | N(t) \ge 1  \}=\frac{\lambda}{\sqrt{\lambda^2+4c^2}} \frac{\sinh \frac{t }{2} \sqrt{\lambda^2+4c^2}}{\sinh \frac{\lambda t}{2}}.
\end{equation}
It is clear that the mean value (\ref{transalpina}) tends to infinity as $t \to \infty$. Furthermore, if $\lambda, c \to \infty$ (so that $\frac{c^2}{\lambda} \to 1$) then $E\{ \cosh \eta_1(t) I_{\{N(t) \ge 1\}}\} \to e^t$. It can also be checked that if $c=0$ then $E\{ \cosh \eta_1(t) I_{\{N(t) \ge 1\}}\} =1$, since the particle never leaves the starting point.

For $k=2$ formula (\ref{arkee}) yields
\begin{eqnarray} \label{karamazov}
E\{ \cosh \eta_2(t)| N(t) \ge 2  \} &=& \frac{\lambda^2 e^{-\frac{\lambda t}{2}}}{c^2 Pr\{ N(t) \ge 2\}} \left( \cosh \frac{t}{2}  \sqrt{\lambda^2+4c^2} \right. \nonumber \\
&-& \left. \frac{\lambda}{\sqrt{\lambda^2+4c^2}} \sinh \frac{t }{2}\sqrt{\lambda^2+4c^2}-e^{-\frac{\lambda t}{2}} \right)\nonumber \\
&=&\frac{\lambda^2}{c^2 Pr\{N(t) \ge 2\}} \left[  E\{ \cosh \eta(t)\} \right.    \\
&-&\left. Pr\{  N(t) \ge 1 \}  E\{ \cosh \eta_1(t) |  N(t)\ge 1 \} -e^{-\lambda t}\right] \nonumber
\end{eqnarray}
Since
\begin{eqnarray}
&&\lim_{c \to 0} \frac{1}{c^2} \left\{ \cosh \frac{t}{2} \sqrt{\lambda^2 +4c^2} - \frac{\lambda}{\sqrt{\lambda^2+4c^2}} \sinh \frac{t}{2} \sqrt{\lambda^2+4c^2} - e^{-\frac{\lambda t}{2}} \right\} \nonumber \\
&&= \frac{e^{\frac{\lambda t}{2}}}{\lambda^2} Pr\{N(t) \ge 2\},
\end{eqnarray}
we have, as expected, that
\begin{equation} \label{13:17 bresaola}
\lim_{c \to 0} E\{ \cosh \eta_2(t) |N(t) \ge 2\}=1.
\end{equation}
Also, when $\lambda \to \infty$, we obtain the same limit as in (\ref{13:17 bresaola}). The expression (\ref{karamazov}) suggests the following decomposition
\begin{eqnarray}
E\{\cosh \eta(t)\} &=&\frac{c^2}{\lambda^2} Pr\{N(t) \ge 2\} E\{ \cosh \eta_2 (t) | N(t) \ge 2 \}\\
&+&Pr\{N(t) \ge 1\} E\{ \cosh \eta_1(t) | N(t) \ge 1 \} +e^{-\lambda t} \nonumber \\
&=&\frac{c^2}{\lambda^2} E\{ \cosh \eta_2 (t) I_{\{ N(t) \ge 2\}} \} + E\{ \cosh \eta_1(t) I_{\{ N(t) \ge 1\} }\} +e^{-\lambda t} \nonumber
\end{eqnarray}
\end{oss}

\begin{oss}
The result in (\ref{9:36coppa}) appears as the mean hyperbolic distance of a motion starting from the origin and running, without returns, until time $t-T_k$, where $T_k$ has a truncated Gamma distribution (Erlang distribution) with density
\begin{equation}
Pr\{T_k \in \mathrm{d}t_k\}= \frac{ \lambda^k e^{-\lambda t_k} t_k^{k-1} }  { (k-1)! Pr\{T_k \le t\} } \mathrm{d}t_k \hspace{1cm} 0<t_k <t.
\end{equation}
In other words we can write (\ref{9:36coppa}) as
\begin{eqnarray} \label{13:53}
E\{ \cosh \eta_k (t) | N(t) \ge k\}&=& E\{E \{ \cosh \eta(t-T_k)\}\} \\
&=& \int_{0}^{t} E\{ \cosh \eta(t-T_k) \} Pr\{T_k \in  \mathrm{d}t_k\}. \nonumber
\end{eqnarray}
Furthermore the expression (\ref{9:36coppa}) contains a fractional integral of order $k$ for the function $g(s)=e^{\lambda s} E(s)$
\begin{eqnarray}
&&E\{ \cosh \eta_k (t) | N(t) \ge k\} \\&&= \frac{\lambda^k}{ \sum_{j=k}^{\infty} \frac{(\lambda t)^j}{j!}} \left\{ \frac{1}{\Gamma(k)} \int_{0}^{t} (t- s)^{k-1} e^{\lambda s} E(s) \mathrm{d} s \right\}. \nonumber
\end{eqnarray}
If the mean value (\ref{13:53}) is taken with respect to
\begin{equation}
Pr\{ T_\nu \in \mathrm{d}s\}=\frac{\lambda^\nu e^{-\lambda s } s^{\nu-1}}{ \Gamma(\nu) Pr\{T_{\nu} \le t\}} \mathrm{d}s\hspace{2cm} 0<s <t,
\end{equation}
then we have that
\begin{eqnarray} \label{14:14}
&&E\{E \{ \cosh \eta(t-T_\nu)\}\}\\
&&=\frac{\lambda^{\nu} e^{-\lambda t}}{Pr\{ T_{\nu} \le t\}} \left\{ \frac{1}{\Gamma(\nu)} \int_{0}^{t} (t-s)^{\nu-1} e^{\lambda s} E(s) \mathrm{d}s \right\} \nonumber
\end{eqnarray}
that also contains a fractional integral of order $\nu$ in the sense of Riemann-Liouville. The expression (\ref{14:14}) can be interpreted as the mean hyperbolic distance at time $t$ where the particle can jump back to the origin at an arbitrary instant (different from the instants of change of direction).
\end{oss}

\section{Motion at Finite Velocity on the Surface of a Three-dimensional Sphere}

Let $P_0$ be a point on the equator of a three-dimensional sphere. Let us assume that the particle starts  moving from $P_0$ along the equator in one of the two possible directions (clockwise or counter-clockwise) with velocity c.

At the first Poisson event (occurring at time $T_1$) it starts moving on the meridian joining the north pole $P_N$ with the position reached at time $T_1$ (denoted by $P_1$) along one of the two possible directions (see Figure \ref{figsfera}).

At the second Poisson event the particle is located at $P_2$ and its distance from the starting point $P_0$ is the length of the hypothenuse of a right spherical triangle with cathetus $P_0 P_1$ and $P_1 P_2$; the hypothenuse belongs to the equatorial circumference through $P_0$ and $P_2$.

Now the particle continues its motion (in one of the two possible directions) along the equatorial circumference orthogonal to the hypothenuse through $P_0$ and $P_2$ until the third Poisson event occurs.

In general, the distance $\mathrm{d}(P_0P_t)$ of the point $P_t$ from the origin $P_0$ is the length of the shortest arc of the equatorial circumference through $P_0$ and $P_t$ and therefore it takes values in the interval $[0,\pi]$. Counter-clockwise motions cover the arcs in $[-\pi, 0]$ so that the distance is also defined in $[0,\pi]$ or in $[-\pi/2, \pi/2]$ with a shift that avoids negative values for the cosine.

\begin{figure}[h]
\centering
\includegraphics[width=7.6cm, height=5cm, clip ]{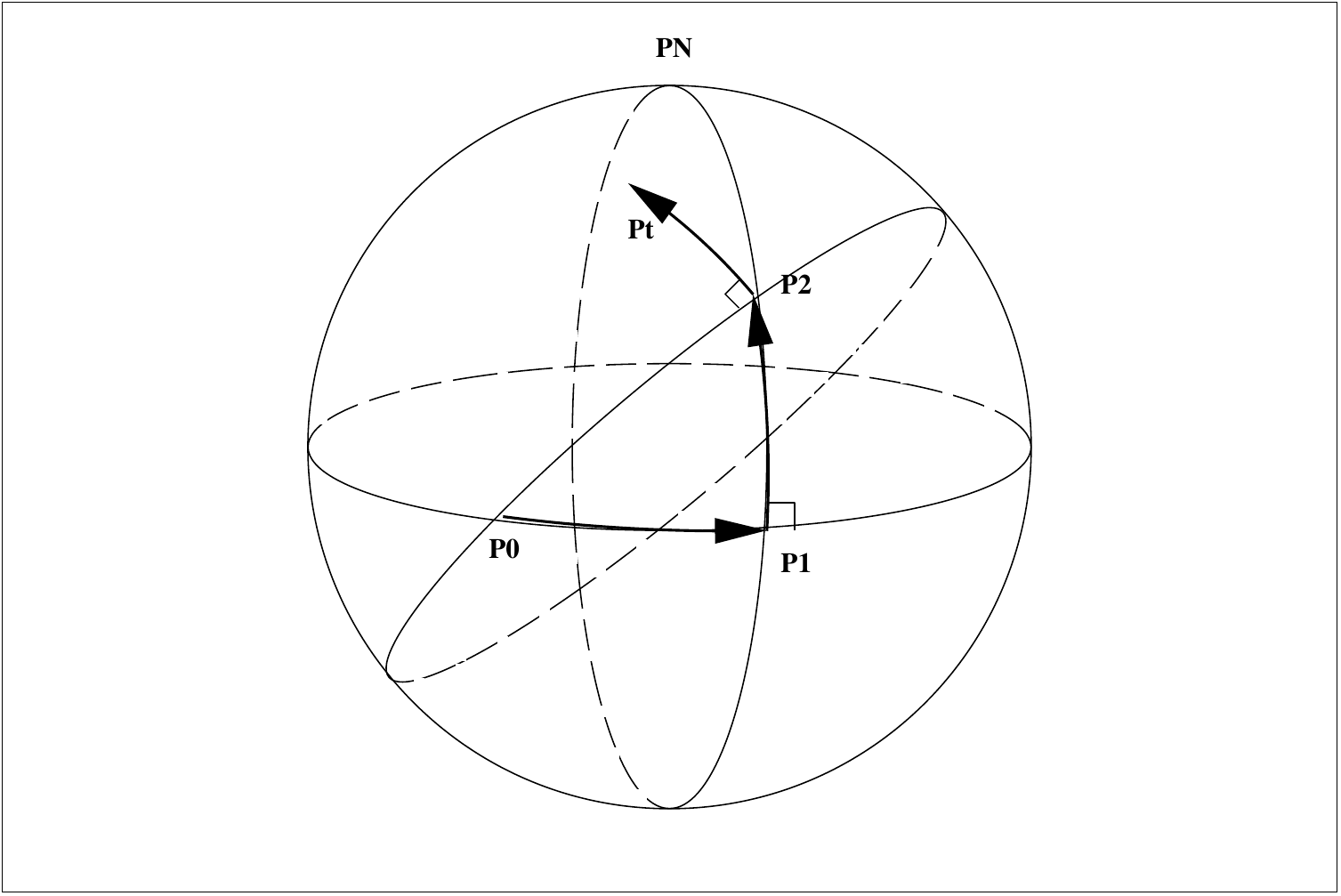}
\caption{Motion on the surface of a three-dimensional sphere. }
\label{figsfera}
\end{figure}

By means of the spherical Pythagorean relationship we have that
the Euclidean distance $\mathrm{d}(P_0 P_2)$ satisfies
\begin{equation}
\cos \mathrm{d}(P_0 P_2)=\cos \mathrm{d}(P_0 P_1) \cos
\mathrm{d}(P_1 P_2)
\end{equation}
and, after three displacements,
\begin{eqnarray}
\cos \mathrm{d}(P_0 P_3)&=&\cos \mathrm{d}(P_0 P_2) \cos \mathrm{d}(P_2 P_3) \nonumber \\
&=&\cos \mathrm{d}(P_0 P_1) \cos \mathrm{d}(P_1 P_2) \cos
\mathrm{d}(P_2 P_3).
\end{eqnarray}
After $n$ displacements the position $P_t$ on the sphere at time $t$ is given by
\begin{equation}
\cos \mathrm{d}(P_0 P_t)=\prod_{k=1}^{n} \cos \mathrm{d}(P_k
P_{k-1}) \cos \mathrm{d}(P_n P_t).
\end{equation}
Since $\mathrm{d}(P_k P_{k-1})$ is represented by the amplitude of
the arc run in the interval $(t_k, t_{k-1})$, it results
$$\mathrm{d}(P_k P_{k-1} )=c(t_k-t_{k-1}).$$

The mean value  $E\{ \cos \mathrm{d}(P_0 P_t)|N(t)=n\}$ is given
by
\begin{eqnarray}
E_n(t)&=&E\{ \cos \mathrm{d}(P_0 P_t)|N(t)=n\}\\
&=&\frac{n!}{t^n}\int_{0}^{t}\mathrm{d}t_1 \int_{t_1}^{t}\mathrm{d}t_2 \cdots \int_{t_{n-1}}^{t}\mathrm{d}t_n \prod_{k=1}^{n+1} \cos c(t_k-t_{k-1}) \nonumber \\
&=&\frac{n!}{t^n} H_n(t), \nonumber
\end{eqnarray}
where $t_0=0$, $t_{n+1}=t$, and
\begin{equation} H_n(t)=\int_{0}^{t}\mathrm{d}t_1 \cdots \int_{t_{n-1}}^{t}\mathrm{d}t_n \prod_{k=1}^{n+1} \cos c(t_k-t_{k-1}).
\end{equation}
The mean value $E\{\cos \mathrm{d}(P_0 P_t)\}$ is given by
\begin{eqnarray}
E(t)&=&E\{\cos \mathrm{d}(P_0 P_t)\}\\
&=&\sum_{n=0}^{\infty} E\{ \cos \mathrm{d}(P_0 P_t)|N(t)=n\} Pr\{N(t)=n\} \nonumber \\
&=& e^{-\lambda t} \sum_{n=0}^{\infty} \lambda^n H_n(t). \nonumber
\end{eqnarray}

By steps similar to those of the hyperbolic case we have that $H_n(t)$, $t\ge 0$, satisfies the difference-differential equation
\begin{equation}
\frac{\mathrm{d}^2}{\mathrm{d}t^2}H_n=\frac{\mathrm{d}}{\mathrm{d}t}H_{n-1}-c^2 H_n,
\end{equation}
where $H_0(t)=\cos ct$, and therefore we can prove the following:
\begin{teo}
The mean value $E(t)=E\{\cos d(P_0 P_t)\} $ satisfies
\begin{equation} \label{fuggono}
\frac{\mathrm{d}^2}{\mathrm{d}t^2}E=- \lambda \frac{\mathrm{d}}{\mathrm{d}t}E-c^2 E
\end{equation}
with initial conditions
\begin{equation} \label{messaggeri}
\left\{
\begin{array}{lr}  E(0) =1, \\
\left. \frac{\mathrm{d}}{\mathrm{d}t}E(t) \right|_{t=0}=0,
\end{array}
\right.
\end{equation}
and has the form
\begin{equation} \label{20:27cena}
E(t)=\left\{
\begin{array}{lr} e^{-\frac{\lambda t}{2}}\left[\cosh \frac{t}{2}\sqrt{\lambda^2-4c^2}+\frac{\lambda}{\sqrt{\lambda^2-4c^2}}\sinh \frac{t}{2}\sqrt{\lambda^2-4c^2} \right] &  0<2c<\lambda,\\
e^{-\frac{\lambda t}{2}}\left[1+\frac{\lambda t}{2} \right] &  \lambda=2c>0,\\
e^{-\frac{\lambda t}{2}}\left[\cos \frac{t}{2}\sqrt{4c^2-\lambda^2}+\frac{\lambda}{\sqrt{4c^2-\lambda^2}} \sin \frac{t}{2}\sqrt{4c^2-\lambda^2} \right] &  2c>\lambda>0.
\end{array}
\right.
\end{equation}
\end{teo}
\Dim
The solution to the problem (\ref{fuggono})-(\ref{messaggeri}) is given by
\begin{eqnarray}
E(t)&=&\frac{e^{- \frac{\lambda t}{2}}}{2}\left[ \left( e^{\frac{t}{2}\sqrt{\lambda^2-4c^2}}+e^{-\frac{t}{2}\sqrt{\lambda^2-4c^2}}\right)\right.
\nonumber \\
&+& \left.  \frac{\lambda}{ \sqrt{\lambda^2-4c^2}}\left( e^{\frac{t}{2}\sqrt{\lambda^2-4c^2}}-e^{-\frac{t}{2}\sqrt{\lambda^2-4c^2}}\right) \right],
\end{eqnarray}
so that (\ref{20:27cena}) emerges.
\Fine

For large values of $\lambda$ the first expression furnishes $E(t)\sim 1$ and therefore the particle hardly leaves the starting point.\\
If $\frac{\lambda}{2}<c$, the mean value exhibits an oscillating behavior; in particular, the oscillations decrease as time goes on, and this means that the particle moves further and further reaching in the limit the poles of the sphere.

\begin{oss}
By assuming that $c$ is replaced by $i c$ in (\ref{20:27cena}) we formally extract from the first and the third expression in  (\ref{20:27cena}) the hyperbolic mean distance (\ref{en}). This is because the space $H_2^+$ can be regarded as a sphere with imaginary radius. Clearly the intermediate case $\lambda=2c$ has no correspondence for the motion on $H_2^+$ because the Poisson rate must be a real positive number.
\end{oss}

\section*{Acknowledgment}
The authors want to thank the referees for their suggestions which improved the first draft.

\end{document}